\documentclass[12pt]{amsart}

\usepackage{amscd,amssymb}
\usepackage{amsthm,amsmath,amssymb}

\usepackage[matrix,arrow]{xy}

\newtheorem{theorem}[equation]{Theorem}

\newtheorem{proposition}[equation]{Proposition}
\newtheorem{lemma}[equation]{Lemma}
\newtheorem{corollary}[equation]{Corollary}
\newtheorem{conjecture}[equation]{Conjecture}

\theoremstyle{definition}
\newtheorem{example}[equation]{Example}
\newtheorem{definition}[equation]{Definition}

\theoremstyle{remark}
\newtheorem{remark}[equation]{Remark}

\author{Ivan Cheltsov}

\title{Birationally superrigid cyclic triple spaces}

\address{\begin{tabbing}\hspace*{28 em}\=\kill
Steklov Institute of Mathematics \>School of Mathematics\\
8 Gubkin street, Moscow 117966   \>University of Edinburgh\\
Russia                           \>Kings Buildings,  Mayfield Road\\
                                 \> Edinburgh EH9 3JZ, UK\\
\texttt{cheltsov@yahoo.com}      \>\\
                                 \>\texttt{I.Cheltsov@ed.ac.uk}
\end{tabbing}}

\thanks{The author is very grateful to M.Grinenko,
V.Is\-kov\-skikh, S.Kudryavtsev, J.Park, Yu.Pro\-kho\-rov,
A.Pukh\-li\-kov and V.Sho\-ku\-rov for fruitful conversations.}

\begin{document}

\begin{abstract}
We prove the birational superrigidity and the nonrationality of a
cyclic triple cover of $\mathbb{P}^{2n}$ branched over a nodal
hypersurface of degree $3n$ for $n\ge 2$. In particular, the
obtained result solves the problem of the birational superrigidity
of smooth cyclic triple spaces. We also consider certain relevant
problems.
\end{abstract}

\maketitle


\section{Introduction.}
\label{sec:1}

The problem of the rationality of an algebraic
variety\footnote{All vareities are assumed to be projective,
normal and defined over $\mathbb{C}$.} is one of the most
interesting problems in algebraic geometry. Global holomorphic
differential forms are natural birational invariants of a smooth
algebraic variety that solve the problem ot the rationality of
algebraic curves and surfaces (see \cite{Za58}, \cite{Is96b}).
However, even in three-dimensional case there are nonrational
varieties that are very close to being rational. In particular,
available discrete invariants does not solve the rationality
problem for higher-dimensional algebraic varieties. For example,
there are nonrational unirational 3-folds (see \cite{IsMa71},
\cite{ClGr72}), which imply that the L\"uroth problem in dimension
$3$ has a negative answer. Unfortunately, there are no known
simple way of proving the nonrationality of higher-dimensional
rationally connected varieties (see \cite{Ko96}, \cite{Is97},
\cite{Ko98}).

There are few known methods of proving the nonrationality of
rationally connected varieties. The finiteness of the group of
birational automorphisms of a smooth quartic 3-fold is proved in
\cite{IsMa71}, which implies its nonrationality. The
nonrationality of a smooth cubic 3-fold is proved in \cite{ClGr72}
through the study of its intermediate Jacobian. The
bi\-ra\-ti\-o\-nal invariance of the torsion subgroup of a group
$H^{3}(\mathbb{Z})$ is used in \cite{ArMu72} to prove the
nonrationality of certain unirational conic bundles. The
nonrationality of a wide class of rationally connected varieties
is proved in \cite{Ko95} by means of the reduction into the
positive characteristic (see \cite{ChWo04}, \cite{Ko96},
\cite{Ko00}).

Every methods of proving the nonrationality of an algebraic
variety has advantages and disadvantages. For example, the method
of the intermediate Jacobian can be applied only to 3-folds, and
except a single csee (see \cite{Ti80a}, \cite{Ti80b}, \cite{We81},
\cite{Ti82}, \cite{Ti86}, \cite{Cl91}) only to 3-folds fibered
into conics (see \cite{Tu72}, \cite{Tu75}, \cite{Be77},
\cite{Tu80}). On the other hand, in the three-dimensional case the
method of the intermediate Jacobian can be often when all other
methods simply can not be used. The degeneration method (see
\cite{Be77}, \cite{Tu80}, \cite{Cl82}, \cite{Bar84}, \cite{Ch04d},
\cite{Ch04e}) shows that sometimes the Griffiths component of the
intermediate Jacobian is the most subtle three-dimensional
birational invariant. For example, an important case of the
rationality criterion of a three-dimensional conic bundle (see
\cite{Is87}, \cite{Is91}, \cite{Is92}, \cite{Is96a}) is proved in
\cite{Sho83} using the intermediate Jacobian method. However,
there are nonrational 3-folds whose group $H^{3}(\mathbb{Z})$ is
trivial (see \cite{Sa82}). In many interesting cases, for example,
for smooth complete intersections, the group $H^{3}(\mathbb{Z})$
has no torsion and, therefore, the method of \cite{ArMu72} can not
be applied (see \cite{CoOj89}, \cite{Pe93}). The method of
\cite{Ko95} works in any dimension, but it proves the
nonrationality of a very general element of an appropriate family.
The technique of \cite{IsMa71} also works in any dimension (see
\cite{Pu00a}), but in general it can be applied only to varieties
that stand too far from the rational ones. For example, it is hard
to believe that one can use the technique of \cite{IsMa71} to get
an example of a smooth deformation of a nonrational variety into a
rational one (see \cite{Tu80}). The latter example is expected to
exist in dimension greater than $3$ (see \cite{Tr84}, \cite{Tr93},
\cite{Has99}, \cite{Has00}).

Let us consider the following notion, which is implicitly
introduced in the paper \cite{IsMa71}, but historically it goes
back to the classical papers \cite{Ne71}, \cite{Fa15},
\cite{Fa47}, but its modern form is considered relatively recently
(see \cite{Co00}, \cite{Pu04a}). Note that the class of terminal
singularities is a higher-dimensional generalization of smooth
points of algebraic surfaces that is closed with respect to the
good birational maps (see \cite{KMM}). The
$\mathbb{Q}$-factoriality simply means that a multiple of every
Weil divisor on a variety is a Cartier divisor. In particular,
every smooth variety has terminal $\mathbb{Q}$-factorial
singularities.

\begin{definition}
\label{definition:superrigidity} A terminal $\mathbb{Q}$-factorial
Fano variety $V$ with $\mathrm{Pic}(V)\cong \mathbb{Z}$ is
bi\-ra\-ti\-o\-nal\-ly superrigid if the following $3$ conditions
hold:
\begin{enumerate}
\item the variety $V$ is not birational to a
fibration\footnote{For every fibration $\tau:Y\to Z$ we assume
that $\mathrm{dim}(Y)>\mathrm{dim}(Z)\ne 0$ and
$\tau_{*}(\mathcal{O}_{Y})=\mathcal{O}_{Z}$.}, whose \hfill\break
generic fiber is a
smooth variety of Kodaira dimension $-\infty$; %
\item the variety $V$ is not birational to a
$\mathbb{Q}$-factorial terminal Fano variety\hfill\break with
Picard
group $\mathbb{Z}$ that is not biregular to $V$; %
\item $\mathrm{Bir}(V)=\mathrm{Aut}(V)$.
\end{enumerate}
\end{definition}

The paper \cite{IsMa71} contains an implicit proof that every
smooth quartic 3-fold in $\mathbb{P}^4$ is birationally superrigid
(see \cite{Co95}). The technique of \cite{IsMa71} can be applied
to certain Fano 3-folds with non-trivial group of birational
automorphisms (see \cite{Is80b}). Therefore one can consider the
following weakened version of the birational superrigidity.

\begin{definition}
\label{definition:rigidity} A terminal $\mathbb{Q}$-factorial Fano
variety $V$ with $\mathrm{Pic}(V)\cong\mathbb{Z}$ is called
bi\-ra\-ti\-o\-nal\-ly rigid if the first two conditions of
Definition~\ref{definition:superrigidity} are satisfied.
\end{definition}

Birationally rigid varieties are nonrational. In particular, there
are no birationally rigid del Pezzo surfaces defined over an
algebraically closed field. However, there are birationally rigid
del Pezzo surfaces over an algebraically non-closed field (see
\cite{Is96b}). Namely, the results of \cite{Ma66} and \cite{Ma67}
imply the birational superrigidity of smooth del Pezzo surfaces of
degree $1$ and the birational rigidity of smooth del Pezzo
surfaces of degree $2$ and $3$ that are defined over a perfect
algebraically non-closed field and have Picard group $\mathbb{Z}$.
In particular, minimal smooth cubic surfaces in $\mathbb{P}^{3}$
are birationally equivalent if and only if they are projectively
equivalent (see \cite{Ma72}).

The birational rigidity and superrigidity can be defined for a
fibration into Fano varieties as well (see \cite{Co00},
\cite{Pu04a}). To be precise, the birational rigidity and
superrigidity can be defined for Mori fibrations (see
\cite{Co95}). Today the birational rigidity is proved for many
smooth 3-folds (see \cite{Is80b}, \cite{Sa80}, \cite{Pu98b},
\cite{Co00}), for many smooth varieties whose dimension is greater
than $3$ (see \cite{Sa82}, \cite{Pu87}, \cite{Pu88a},
\cite{Pu98a}, \cite{Ch00b}, \cite{Pu00b}, \cite{Pu00c},
\cite{Pu00d}, \cite{Pu01}, \cite{Sob02}, \cite{Pu02a},
\cite{Pu02b}, \cite{Pu02c}, \cite{dFEM03}, \cite{Ch03b},
\cite{Pu03b}), and for many singular varieties (see \cite{Pu88b},
\cite{Pu97}, \cite{CPR}, \cite{Co00}, \cite{Me03}, \cite{Pu03a},
\cite{ChPa04}, \cite{Ch04h}). For some birationally nonrigid
algebraic varieties it is possible to find all Mori fibrations
birational to them (see \cite{Gr98}, \cite{CoMe02}, \cite{Gr03},
\cite{Gr04}). Unfortunately, despite the obvious success in this
area of algebraic geometry there are many still unsolved relevant
classical problems such as finding the generators of the group
$\mathrm{Bir}(\mathbb{P}^{3})$ or finding the generators of the
group of birational automorphisms of a smooth cubic 3-fold. The
solution of the latter problem is announced in the classical paper
\cite{Fa47}, but the proof contains many gaps.

In the given paper we will prove the following result.

\begin{theorem}
\label{theorem:main}%
Let $\pi:X\to\mathbb{P}^{2n}$ be a cyclic triple
cover\,\footnote{A finite morphism of degree $3$ that induces the
cyclic extension of the fields of rational functions.} such that
$\pi$ is branched over a hypersurface $S\subset\mathbb{P}^{2n}$ of
degree $3n$, $n\ge 2$ and the hypersurface $S$ has at most
ordinary double points. Then $X$ is a terminal
$\mathbb{Q}$-factorial Fano variety with
$\mathrm{Pic}(X)\cong\mathbb{Z}$ such that $X$ is birationally
superrigid, the group $\mathrm{Bir}(X)$ is finite and for
sufficiently  general hypersurface $S\subset\mathbb{P}^{2n}$ it is
isomorphic to $\mathbb{Z}_{3}$. In particular, the variety $X$ is
nonrational.
\end{theorem}

\begin{remark}
\label{remark:hypersurface-in-wps}%
Under the conditions of Theorem~\ref{theorem:main}, the variety
$X$ can be considered as a hypersurface in the weighted projective
space $\mathbb{P}(1^{2n+1}, n)$ of degree $3n$ given by the
equation
$$
y^{3}=f_{3n}(x_{0},\ldots,x_{2n})\subset\mathbb{P}(1^{2n+1},n)\cong\mathrm{Proj}(\mathbb{C}[x_{0},\ldots,x_{2n},y]),%
$$
where $f_{3n}$ is a homogenesous polynomial of degree $3n$ (see
\cite{Mi85}, \cite{Fuj88}, \cite{Tan91}, \cite{Tan92},
\cite{Tan02}), and $\pi:X\to\mathbb{P}^{2n}$ is a restriction of
the natural projection
$\mathbb{P}(1^{2n+1},n)\dashrightarrow\mathbb{P}^{2n}$ induced by
the embedding of the graded algebras
$\mathbb{C}[x_{0},\ldots,x_{2n}]\subset\mathbb{C}[x_{0},\ldots,x_{2n},y]$.
Moreover, the hypersurface $S\subset\mathbb{P}^{2n}$ is given by
the equation $f_{3n}(x_{0},\ldots,x_{2n})=0$.
\end{remark}

\begin{remark}
\label{remark:non-rigid-triple-covers}%
Consider a cyclic triple cover $\pi:X\to\mathbb{P}^{k}$ such that
$\pi$ is branched over a nodal hypersurface
$S\subset\mathbb{P}^{k}$ of degree $3n$ and $k\ge 3$. Then $X$ is
not birationally superrigid in the case when $k<2n$, because it
has pencils of varieties of Kodaira dimension $-\infty$. On the
other hand, the Kodaira dimension of the variety $X$ is
non-negative when $k>2n$ and the variety $X$ is not even uniruled
in this case. Therefore, all birationally superrigid smooth cyclic
triple covers are described by Theorem~\ref{theorem:main}.
\end{remark}

\begin{corollary}
\label{corollary:fields}%
Let $f(x_{0},\ldots,x_{2n})$ be a homogeneous polynomial of degree
$3n$ such that
$$
f(x_{0},\ldots,x_{2n})=0\subset\mathbb{P}^{2n}\cong\mathrm{Proj}(\mathbb{C}[x_{0},\ldots,x_{2n}])%
$$
is a nodal or smooth hypersurface. Then
$\mathbb{C}(\nu_{1},\ldots,\nu_{2n})\sqrt[3]{f(1,\nu_{1},\ldots,\nu_{2n})}$
is a purely transcendental extension of the field $\mathbb{C}$ if
and only if the equality $n=1$ holds.
\end{corollary}

\begin{example}
\label{example:nodal-I} Let $X$ be a hypersurface in
$\mathbb{P}(1^{2n+1}, n)$ of degree $3n$ whose equation is
$$
y^{3}=\sum_{i=0}^{2n}x_{i}^{3n}\subset\mathbb{P}(1^{2n+1},n)\cong\mathrm{Proj}(\mathbb{C}[x_{0},\ldots,x_{2n},y]),%
$$
and $n\ge 2$. Then the projection
$\pi:X\to\mathbb{P}^{2n}\cong\mathrm{Proj}(\mathbb{C}[x_{0},\ldots,x_{2n}])$
is a cyclic triple cover branched over a smooth hypersurface
$\sum_{i=0}^{2n}x_{i}^{3n}=0$, the variety $X$ is birationally
superrigid by Theorem~\ref{theorem:main} and
$$
\mathrm{Bir}(X)=\mathrm{Aut}(X)\cong\mathbb{Z}_{3}\oplus\mathrm{Aut}(\sum_{i=0}^{2n}x_{i}^{3n}=0)\cong\mathbb{Z}_{3}\oplus(\mathbb{Z}_{3n}^{2n}\rtimes\mathrm{S}_{2n+1}),%
$$
where $\mathrm{S}_{2n+1}$ is a symmetric group (see \cite{We79},
\cite{Shi82}, \cite{Shi88}, \cite{Kon02}). Hence $X$ is
nonrational and
$\mathbb{C}(\nu_{1},\ldots,\nu_{2n})\sqrt[3]{1+\sum_{i=1}^{2n}\nu_{i}^{3n}}$
is not a purely transcendental extension of $\mathbb{C}$.
\end{example}

\begin{example}
\label{example:nodal-II} Let $X$ be a hypersurface
$\mathbb{P}(1^{2n+1}, n)$ of degree $3n$ whose equation is
$$
y^{3}=\sum_{i=1}^{n}a_{i}(x_{0},\ldots,x_{2n})x_{i}\subset\mathbb{P}(1^{2n+1},n)\cong\mathrm{Proj}(\mathbb{C}[x_{0},\ldots,x_{2n},y]),%
$$
where $a_{i}$ is a sufficiently general homogeneous polynomial of
degree $3n-1$. Then the natural projection
$\pi:X\to\mathbb{P}^{2n}$ is a cyclic triple cover such that $\pi$
is branched over a nodal hypersurface $S\subset\mathbb{P}^{2n}$ of
degree $3n$, which is given by the equation
$$
\sum_{i=1}^{n}a_{i}x_{i}=0\subset\mathbb{P}^{2n}\cong\mathrm{Proj}(\mathbb{C}[x_{0},\ldots,x_{2n}])
$$
and which has $(3n-1)^{n}$ ordinary double points. The variety $X$
is birationally superrigid and nonrational for $n\ge 2$ by
Theorem~\ref{theorem:main}, and the group $\mathrm{Bir}(X)$ is
finite.
\end{example}

\begin{example}
\label{example:nodal-III} Let $X$ be a hypersurface in
$\mathbb{P}(1^{2n+1}, n)$ of degree $3n$ whose equation is
$$
y^{3}=\sum_{i=1}^{n}a_{i}(x_{0},\ldots,x_{2n})b_{i}(x_{0},\ldots,x_{2n})\subset\mathbb{P}(1^{2n+1},n)\cong\mathrm{Proj}(\mathbb{C}[x_{0},\ldots,x_{2n},y]),%
$$
where $a_{i}$ and $b_{i}$ are sufficiently general homogeneous
polynomials of degree $2n$ and $n$ respectively. Then the natural
projection $\pi:X\to\mathbb{P}^{2n}$ is a cyclic triple cover
branched over a nodal hypersurface $S\subset\mathbb{P}^{2n}$ of
degree $3n$ having $2^{n}n^{2n}$ ordinary double points, the
variety $X$ is birationally superrigid and nonrational for $n\ge
2$ by Theorem~\ref{theorem:main}, and the group $\mathrm{Bir}(X)$
is finite.
\end{example}

The main reason why the variety $X$ in Theorem~\ref{theorem:main}
is birationally superrigid is the following: the anticanonical
degree $(-K_{X})^{\mathrm{dim}(X)}=3$ of the variety $X$ is very
small and the singularities of the variety $X$ are relatively
mild. Roughly speaking, a Fano variety must become \emph{more
rational} when the anticanonical degree getting bigger and the
singularities getting worse. This general principle may not
necessary be true in certain extremely singular cases (see
\cite{Ch97b}). However, it follows from the classification that a
smooth Fano 3-fold is rational if its degree is bigger than $24$
(see \cite{IsPr99}). Singular Fano 3-folds are not classified even
in the case when their anticanonical divisors are Cartier divisors
(see \cite{Ch99}, \cite{Pr04a}, \cite{JaRa04}), but many examples
affirm the intuition in the singular case as well (see \cite{CPR},
\cite{CoMe02}, \cite{Ch97a}, \cite{Ch04b}, \cite{Ch04e},
\cite{ChPa04}). Therefore the nonrationality of the variety $X$ in
Theorem~\ref{theorem:main} is very natural.

Due to natural reasons, it makes sense to consider birational
superrigidity and birational rigidity only for Mori fibrations
(see \cite{Co95}). In particular, in the case of Fano varieties we
must assume that for a given Fano variety its singularities are
$\mathbb{Q}$-factorial and its rank of the Picard group is $1$.
Many examples suggest that a Fano variety may not be birationally
rigid if its degree is not sufficiently small. Moreover, it is
intuitively clear the quantitive characteristics of singularities
(number of isolated singular points or anticanonical degree of the
corresponding subvarieties of singular points) is important only
to provide the $\mathbb{Q}$-factoriality condition (see
\cite{CoMe02}, \cite{Me03}, \cite{ChPa04}, \cite{Ch04e},
\cite{Ch04h}). On the other hand, the qualitative characteristics
of singularities (multiplicity and analytical local type) can
crucially influence the birational geometry of a Fano variety (see
\cite{Co00}, \cite{CoMe02}).

Unfortunately, all existent proofs of the birational rigidity or
birational superrigidity of a Fano variety crucially  depend on
the projective geometry of the given variety related to the
anticanonical map. It is natural to expect that some claims on
birational rigidity can be proven without implicit usage of the
properties of the anticanonical ring. For example, we expect that
the following is true (cf. \cite{Pu02a}, \cite{dFEM03}).

\begin{conjecture}
\label{conjecture:birational-rigidity}%
Let $X$ be a smooth Fan variety of dimension $k$ such that
$\mathrm{Pic}(X)\cong\mathbb{Z}$ and $(-K_{X})^{k}\le 2(k-1)$.
Then $X$ is birationally rigid.
\end{conjecture}

It should be pointed out that
Conjecture~\ref{conjecture:birational-rigidity} is proved only in
dimension $3$ through the explicit classification of smooth Fano
3-folds (see \cite{IsPr99}). It is very possible that the proof of
Conjecture~\ref{conjecture:birational-rigidity} can be extremely
hard. On the other hand, it is very natural to expect that the
following weakened version of the
Conjecture~\ref{conjecture:birational-rigidity} can be proved
relatively soon using methods of \cite{Co00}, \cite{EinLa93},
\cite{Ka97}.

\begin{conjecture}
\label{conjecture:birational-rigidity-simple}%
Let $X$ be a smooth Fano variety of dimension $k$ such that
$\mathrm{Pic}(X)\cong\mathbb{Z}$ and $(-K_{X})^{k}=1$. Then $X$ is
birationally superrigid.
\end{conjecture}

\begin{remark}
\label{remark:non-closed-fields}%
It is well known that any statement on birational rigidity remains
true over any field of definition of the considered varieties with
a single exception. Namely, the characteristic of the field of
definition must be zero in order to use the Kawamata--Viehweg
vanishing theorem (see \cite{Ka82}, \cite{Vi82}). However, in the
case of algebraic surfaces it is enough to assume that the field
of definition is just perfect (see \cite{Ma66}, \cite{Ma67}).
Moreover, one can consider equivariant version of any statement on
birational rigidity when the acting group if finite (see
\cite{Is79}, \cite{Gi81}, \cite{Is96b}). The latter can be used in
classification of all nonconjugate finite subgroups of
corresponding groups of birational automorphisms (see
\cite{Is03}).
\end{remark}

It should be pointed out that the nonrationality and the
non-ruledness of a cyclic triple cover of $\mathbb{P}^{2n}$
branched over a very general\footnote{A complement to a countable
union of Zariski dense subsets in moduli.} smooth hypersurface of
degree $3n$ with $n\ge 2$ are implied by Theorem~5.13 in
\cite{Ko96} that claims the following.

\begin{theorem}
\label{theorem:Janos}%
Let $\xi:V\to\mathbb{P}^{k}$ be a cyclic cover of prime degree
$p\ge 2$ branched over a very general hypersurface
$F\subset\mathbb{P}^{k}$ of degree $pd$ such that $k\ge 3$ and
$d>{\frac {k+1}{p}}$. Then $V$ is nonruled and, in particular, the
variety $V$ is nonrational.
\end{theorem}

In the conditions and notations of Theorem~\ref{theorem:main}, it
is natural to ask how many singular points can $X$ have. The
singular points of the variety $X$ are in one-to-one
correspondence with oddinary double points of the hypersurface
$S\subset\mathbb{P}^{2n}$ of degree $3n$. Therefore, the best
known bound is due to \cite{Va83}. Namely, the number of singular
points of $X$ does not exceed the Arnold number
$\mathrm{A}_{2n}(3n)$, where $\mathrm{A}_{2n}(3n)$ is a number of
points $(a_{1},\ldots,a_{2n})\subset\mathbb{Z}^{2n}$ such that the
inequalities
$$
3n^{2}-3n+2\le \sum_{i=1}^{2n}a_{i}\le 3n^{2}
$$
hold and $a_{i}\in(0,3n)$. In particular, the number of singular
points of the variety $X$ does not exceed $320$, $115788$ and
$85578174$ when $n=2$, $3$ and $4$ respectively. However, this
bound seems not to be sharp for $n\gg 0$ (see \cite{Be79a},
\cite{St78}, \cite{CaCe82}, \cite{vSt93}, \cite{Ba96},
\cite{JaRu97}, \cite{Wa98}).

\begin{remark}
\label{remark:rationally-connected} It is well known that the
variety $X$ in Theorem~\ref{theorem:main} is a rationally
connected variety (see \cite{KoMiMo1}, \cite{KoMiMo3},
\cite{KoMiMo3}, \cite{Ko96}). Namely, there is an irreducible
rational curve on the variety $X$ passing through any two
sufficiently general points of $X$.
\end{remark}

The geometrical meaning of Theorem~\ref{theorem:main} has the same
nature as the Noether theorem that claims that the group
$\mathrm{Bir}(X)$ is generated by the Cremona involution and
projective automorphisms (see \cite{Ne71}, \cite{Is80b},
\cite{Co95}). The Noether theorem is related to many interesting
problems. For example, the Noether theorem is related to the
problem of birational classification of plane elliptic pencils.
Originally it was considered in \cite{Ber77}, but later the ideas
of \cite{Ber77} were put into proper and correct form in the paper
\cite{Do66} that proves that any plane elliptic pencil can be
birationally transformed into a special plane elliptic pencil,
so-called Halphen pencil (see \S 5.6 in \cite{CosDo89}), which was
studied in  \cite{Hal82}. A similar problem can be considered for
the variety $X$ in Theorem~\ref{theorem:main}. Namely, we prove
the following result.

\begin{theorem}
\label{theorem:second} Under the conditions of
Theorem~\ref{theorem:main}, the variety $X$ is not birational to
any elliptic fibration.
\end{theorem}

Birational transformations into elliptic fibrations were used in
\cite{BoTsch99}, \cite{BoTsch00}, \cite{HaTsch00} in the proof of
the potential density\footnote{The set of rational points of a
variety $V$ defined over a number field $\mathbb{F}$ is called
potentially dense if for a finite extension of fields
$\mathbb{K}\slash\mathbb{F}$ the set of $\mathbb{K}$-rational
points of the variety $V$ is Zariski dense.} of rational points on
smooth Fano 3-folds, where the following result was proved.

\begin{theorem}
\label{theorem:rational-density}%
Rational points are potentially dense on all smooth Fano 3-folds
with a possible exception of a double cover of $\mathbb{P}^{3}$
ramified in a smooth sextic surface.
\end{theorem}

The existence of a possible exception in
Theorem~\ref{theorem:rational-density} is explained by the
following result proved in \cite{Ch00a}: a smooth sextic double
solid is the only smooth Fano 3-fold that is not birationally
isomorphic to an elliptic fibration (see \cite{IsPr99}). It should
be pointed out that a double cover of $\mathbb{P}^{3}$ branched
over a sextic having one ordinary double point can be birationally
transformed into an elliptic fibration in a unique way (see
\cite{Ch01a}) and rational points on such 3-fold are potentially
dense (see \cite{ChPa04}).

\begin{remark}
\label{remark:counterexample}%
Let $\pi:X\to\mathbb{P}^{4}$ be a cyclic triple cover such that
$\pi$ is branched over a hypersurface $S\subset\mathbb{P}^{4}$ of
degree $6$, $n\ge 2$, and $S$ has one ordinary singular point
$O\in S$ of multiplicity $3$. Then the projection
$\gamma:\mathbb{P}^{4}\dasharrow\mathbb{P}^{3}$ from $O$ induces a
rational map $\gamma\circ\pi$ such that the normalization of the
generic fiber of $\gamma\circ\pi$ is an elliptic curve. In
particular, the variety $X$ does not satisfy the conditions of
Theorem~\ref{theorem:main}. Namely, $S$ is not nodal.
\end{remark}

The nodality condition in Theorems~\ref{theorem:main} and
\ref{theorem:second} is rather natural. Indeed, ordinary double
points are the simplest singularities of algebraic varieties and
the geometry of nodal varieties is related to many interesting
problems (see \cite{To36}, \cite{Cl83}, \cite{Fi87}, \cite{We87},
\cite{JSV90}, \cite{Bo90}, \cite{Pet98}, \cite{Cy99}, \cite{En99},
\cite{Cy01}). However, we can consider a wider class of
singularities in the problems similar to the claim of
Theorems~\ref{theorem:main}. The proofs of
Theorems~\ref{theorem:main} and \ref{theorem:second} together with
the inequality for global log canonical thresholds (see
\cite{Ch01b}, \cite{ChPa02}, \cite{EM01}) give a proof of the
following simple generalization of Theorems~\ref{theorem:main} and
\ref{theorem:second}.

\begin{theorem}
\label{theorem:third}%
Let $\pi:X\to\mathbb{P}^{2n}$ be a cyclic triple cover such that
$\pi$ is branched over a hypersurface $S\subset\mathbb{P}^{2n}$ of
degree $3n$, $n\ge 2$ and the only singularities of $S$ are
ordinary double and triple points. Namely, the multiplicity of any
singular point of $S$ does not exceed $3$ and the projectivization
of the tangent cone to the hypersurface $S$ at this point is
smooth. Then $X$ is a Fano variety with $\mathbb{Q}$-factorial
terminal singularities, $\mathrm{Pic}(X)\cong\mathbb{Z}$, the
variety $X$ is birationally superrigid, and the group
$\mathrm{Bir}(X)$ is finite. Moreover, the only way to
birationally transform $X$ into an elliptic fibration is by means
of the construction in Remark~\ref{remark:counterexample}, which
implies $n=2$ and $S$ has a triple point.
\end{theorem}

Therefore, it follows from Theorem~\ref{theorem:third} that the
methods of \cite{BoTsch99}, \cite{BoTsch00}, \cite{HaTsch00} can
not be used to prove the potential density of rational points on
the variety $X$ in Theorem~\ref{theorem:third} in the case when
the variety $X$ is defined over a number field, with a single
exception of a cyclic triple cover of $\mathbb{P}^{4}$ branched
over a hypersurface of degree $6$ having at least one triple
point. It should be pointed out that the geometrical
unirationality of a variety defined over a number field implies
the potential density of rational points. Therefore, if rational
points are not potentially dense on some of the considered cyclic
triple covers, then it is rationally connected but not
unirational! On the other hand, as of today there is no known
example of a rationally connected variety that is not unirational
(cf. Conjecture 4.1.6 in \cite{Ko98}). Therefore, it is natural to
expect that the methods of \cite{BoTsch99}, \cite{BoTsch00} and
\cite{HaTsch00} can be applied to prove potential density of
rational points of a cyclic triple cover of $\mathbb{P}^{4}$ which
is defined over a number field and branched over a hypersurface of
degree $6$ having at least one singular point of multiplicity $3$.
We will prove the this statement in the general case only. Namely,
we will prove the following result using the  method of
\cite{BoTsch99}, \cite{BoTsch00}, \cite{HaTsch00}.

\begin{theorem}
\label{theorem:five}%
Let $\pi:X\to\mathbb{P}^{4}$ be a cyclic triple cover branched
over a sufficiently general\,\footnote{A complement to a Zariski
closed subset in moduli.} hypersurface $S\subset\mathbb{P}^{4}$ of
degree $6$ such that $S$ is defined over a number field and $S$
has an ordinary triple point. Then rational points are
potentially dense on $X$.
\end{theorem}

Actually, our methods can be used to prove the following result.
Let us remind that canonical singularities are higher-dimensional
generalization of Du Val singularities of algebraic surfaces (see
\cite{KMM}).

\begin{theorem}
\label{theorem:forth}%
Under the conditions of Theorem~\ref{theorem:main} or
Theorem~\ref{theorem:third}, let $\rho:X\dasharrow V$ be a
birational map such that $V$ is a Fano variety with canonical
singularities. Then $\rho$ is an isomorphism.
\end{theorem}

The claim of Theorem~\ref{theorem:forth} is a generalization of
one of the claims of Theorem~\ref{theorem:main}. However, we think
that Theorem~\ref{theorem:forth} has certain importance. For
example, the similar claim for smooth minimal cubic surfaces
defined over an algebraically non-closed field (see \cite{Ch00a})
generalizes the classical birational classification (see
\cite{Ma72}) in the following way: a smooth minimal cubic surface
in $\mathbb{P}^{3}$ is birational to a cubic surface in
$\mathbb{P}^{3}$ with Du Val singularities if and only if they are
projectively equivalent. Moreover, the expanded version of the
latter claim (see \cite{Ch00a}) gives a description of all finite
subgroups of the group of birational automorphisms of a smooth
minimal cubic surface (see \cite{Ch04g}), which answers
Question~1.10 in the book \cite{Ma72}. The latter problem was
originally solved in \cite{Ka77} by group-theoretic methods using
the the explicit description of the group of birational
automorphisms of a smooth minimal cubic surface obtained in
\cite{Ma67} and \cite{Ma72}.

\begin{remark}
\label{remark:remark} The claims similar to
Theorems~\ref{theorem:second} and \ref{theorem:forth} are proved
for many algebraic varieties (see \cite{Ch00a}, \cite{Ch00b},
\cite{Ch00c}, \cite{Ch01a}, \cite{Ry02}, \cite{Ch03a},
\cite{Ch03b}, \cite{Ch04a}, \cite{Ch04c}, \cite{Ch04f},
\cite{ChPa04}).
\end{remark}

Double covers of projective spaces are generalizations of
hyperelliptic curves, triple covers of projective spaces are
generalizations of trigonal curves. However, triple covers are not
necessary Galois covers. The study of discrete invariants of
cyclic covers of $\mathbb{P}^{2}$ goes back to \cite{Com11},
\cite{Za29}, \cite{Za31}, which was continued in the papers
\cite{Iv70}, \cite{Lib82}, \cite{Sak94}, \cite{To99},
\cite{CatCi93} and \cite{Ku95}. Certain questions related to
triple covers of algebraic surfaces were considered in
\cite{To91}, \cite{Tan91}, \cite{Tan92}. The topological questions
related to covers of projective spaces were considered in
\cite{La80} and \cite{FuLa81}. Structural results related to
triple covers were obtained in \cite{Mi85}, \cite{Fuj88},
\cite{Par89}, \cite{CasEck96}, \cite{Tan01}, \cite{Tan02},
\cite{FaeSti02}. Some results of sporadic nature were obtained in
\cite{vdWa78}, \cite{Par91}, \cite{Man97}. In the framework of
birational geometry triple covers of projective spaces were
considered in \cite{LanLi90} and \cite{LanPaSo94}. The
nonrationality of general cyclic covers of projective spaces were
considered in \cite{Ko96} (see Theorem~\ref{theorem:Janos}).

\section{Movable log pairs.}
\label{sec:2}

In this section we will consider properties of so-called movable
log pairs that were introduced in \cite{Al91}. Movable log pair
were used implicitly in \cite{Ne71}, \cite{Fa15}, \cite{Fa47},
\cite{IsMa71}.

\begin{definition}
\label{definition:movable-log-pair} A movable log pair $(X,
M_{X})$ is pair consisting of a variety $X$ and a movable boundary
$M_{X}$, where $M_{X}=\sum_{i=1}^{n} a_{i}\mathcal{M}_{i}$ is a
formal finite linear combination of linear systems
$\mathcal{M}_{i}$ on variety $X$ such that the base locus of every
$\mathcal{M}_{i}$ has codimension at least $2$ in $X$ and
$a_{i}\in \mathrm{Q}_{\ge 0}$.
\end{definition}

It is clear that every movable log pair can be considered as a
usual log pair with an effective boundary whose components does
not have multiplicities greater than $1$ by replacing every linear
system either by its general element or by the appropriate
weighted sum of its general elements. In particular, for a given
movable log pair $(X, M_{X})$ we may consider movable boundary
$M_{X}$ as an effective divisor. Thus the numerical intersection
of the movable boundary $M_{X}$ with curves on the variety $X$ is
well defined in the case when the variety $X$ is
$\mathbb{Q}$-factorial. Hence we can consider the formal sum
$K_{X}+M_{X}$ as a log canonical divisor of the movable log pair
$(X, M_{X})$. In the rest of this section we assume that all log
canonical divisors are $\mathbb{Q}$-Cartier divisors.

\begin{remark}
\label{remark:square} For a movable log pair $(X, M_{X})$ the
self-intersection $M_{X}^{2}$ can be considered as a well-defined
effective codimension-two cycle in the case when the singularities
of the variety $X$ are $\mathbb{Q}$-factorial.
\end{remark}

The image of a movable boundary under a birational map is
naturally well defined, because base loci of the components of a
movable boundary do not contain divisors.

\begin{definition}
\label{definition:birational-movable-log-pairs} Movable log pairs
$(X, M_{X})$ and $(Y, M_{Y})$ are called birationally equivalent
if there is a birational map $\rho:X\dasharrow Y$ such that
$M_{Y}=\rho(M_{X})$.
\end{definition}

The standard notions such as discrepancies, terminality,
cano\-ni\-city, log terminality and log cano\-ni\-city can be
defined for movable log pairs in a similar way as they are defined
for usual log pairs (see \cite{KMM}).

\begin{definition}
\label{definition:canonical-singularities} A movable log pair $(X,
M_{X})$ has canonical (terminal respectively) singularities if for
every birational morphism $f:W\to X$ there is an equivalence
$$
K_{W}+f^{-1}(M_{X})\sim_{\mathbb{Q}} f^{*}(K_{X}+M_{X})+\sum_{i=1}^{n}a(X, M_{X}, E_{i})E_{i}%
$$
such that every rational number $a(X, M_{X}, E_{i})$ is
non-negative (positive respectively), where $E_{i}$ is an
$f$-exceptional divisor. The rational number $a(X, B_{X}, E_{i})$
is called a discrepancy of the movable log pair $(X, B_X)$ in the
$f$-exceptional divisor $E_i$.
\end{definition}

\begin{example}
\label{example:reduced-boundary} Let $X$ be a 3-fold and
$\mathcal{M}$ be a linear system on $X$ such that the base locus
of the linear system $\mathcal{M}$ has codimension at least $2$.
Then the log pair $(X, \mathcal{M})$ has terminal singularities if
and only if the linear system $\mathcal{M}$ has only isolated
simple base points, which are smooth points of the 3-fold $X$.
\end{example}

\begin{remark}
\label{remark:LMMP} The application of Log Minimal Model Program
(see \cite{KMM}) to a movable log pair having canonical or
terminal singularities preserves the canonicity or terminality
respectively.
\end{remark}

Singularities of a movable log pair coincide with the
singularities of the variety outside of the base loci of the
components of the movable boundary. Therefore the existence of a
resolution of singularities (see \cite{Hi64}) implies that every
movable log pair is birationally equivalent to a log pair with
canonical or terminal singularities.

\begin{definition}
\label{definition:center} A proper irreducible subvariety
$Y\subset X$ is called a center of canonical singularities of a
movable log pair $(X, M_{X})$ if there is a birational morphism
$f:W\to X$ and an $f$-ex\-cep\-tional divisor $E_{1}\subset W$
such that
$$
K_{W}+f^{-1}(M_{X})\sim_{\mathbb{Q}} f^{*}(K_{X}+M_{X})+\sum_{i=1}^{k} a(X, M_{X}, E_{i})E_{i},%
$$
where $a(X, M_{X}, E_{i})\in \mathbb{Q}$, $E_{i}$ is an
$f$-exceptional divisor, $a(X, M_{X}, E_{1})\le 0$, $f(E_{1})=Y$.
\end{definition}

\begin{definition}
\label{definition:set-of-centers} The set $\mathbb{CS}(X, M_{X})$
is a set of all centers of canonical sin\-gu\-la\-ri\-ties of a
movable log pair $(X, M_{X})$ and $CS(X, M_{X})$ is a
set-theoretic union of all centers of canonical singularities of
the movable log pair $(X, M_{X})$.
\end{definition}

In particular, a log pair $(X, M_{X})$ is terminal if and only if
$\mathbb{CS}(X, M_{X})=\emptyset$.

\begin{remark}
\label{remark:flexability-terminality} Let $(X, M_{X})$ be a log
pair with terminal singularities. Then the sin\-gu\-la\-rities of
the log pair $(X, \epsilon M_{X})$ are terminal for any small
enough rational number $\epsilon>1$.
\end{remark}

\begin{remark}
\label{remark:canonical-centers-of-codimension-two} Let $(X,
M_{X})$ be a movable log pair and $Z\subset X$ be a proper
irreducible subvariety such that $X$ is smooth at the generic
point of the subvariety $Z$. Then elementary properties of blow
ups imply
$$
Z\in \mathbb{CS}(X, M_{X})\Rightarrow \mathrm{mult}_{Z}(M_{X})\geq 1%
$$
and in the case $\mathrm{codim}(Z\subset X)=2$ we have
$\mathrm{mult}_{Z}(M_{X})\geq 1\Rightarrow Z\in \mathbb{CS}(X,
M_{X})$.
\end{remark}

\begin{remark}
\label{remark:canonical-reduction} Let $(X, M_{X})$ be a movable
log pair, $H$ be a general hyperplane section of the variety $X$,
and $Z\in \mathbb{CS}(X, M_{X})$ such that $\mathrm{dim}(Z)\ge 1$.
Then $Z\cap H\in \mathbb{CS}(H, M_{X}\vert_{H})$.
\end{remark}

\begin{definition}
\label{definition:Kodaira-dimension} For a movable log pair $(X,
M_{X})$ consider any birationally equivalent movable log pair $(W,
M_{W})$ such that its singularities are canonical. Let $m$ be a
natural number such that $m(K_{W}+M_{W})$ is a Cartier divisor.
The Kodaira dimension $\kappa(X, M_{X})$ of the log pair $(X,
M_{X})$ is the maximal dimension of the image
$\phi_{|nm(K_{W}+M_{W})|}(W)$ for $n\gg 0$ in the case when
$|n(K_{W}+M_{W})|$ is not empty for some $n$. In the case when
complete linear systems $|n(K_{W}+M_{W})|$ are empty for all $n\gg
0$ we simply put $\kappa(X, M_{X})=-\infty$.
\end{definition}

\begin{lemma}
\label{lemma:Kodaira-dimension-is-OK} The Kodaira dimension of a
movable log pair is well-defined, namely, it does not depend on
the choice of the birationally equivalent movable log pair having
canonical singularities in
Definition~\ref{definition:Kodaira-dimension}.
\end{lemma}

\begin{proof}
Let $(X, M_{X})$ and $(Y, M_{Y})$ be movable log pairs having
canonical singularities such that we have $M_{X}=\rho(M_{Y})$ for
some birational map $\rho:Y\dasharrow X$. Take positive integer
$m$ such that $m(K_{X}+M_{X})$ and $m(K_{Y}+M_{Y})$ are Cartier
divisors. To conclude the claim it is enough to show that
$\phi_{|nm(K_{X}+M_{X})|}(X)=\phi_{|nm(K_{Y}+M_{Y})|}(Y)$ for
$n\gg 0$ or
\[|nm(K_{X}+M_{X})|=|nm(K_{Y}+M_{Y})|=\emptyset \mbox{ for all } n\in\mathbb{N}.\]

Let us consider a commutative diagram
\[\xymatrix{
&&W\ar@{->}[ld]_{g}\ar@{->}[rd]^{f}&&\\%
&X\ar@{-->}[rr]_{\rho}&&Y&}\] %
such that $W$ is smooth, $g:W\to X$ and $f:W\to Y$ are birational
morphisms. Then
$$
K_{W}+M_{W}\sim_{\mathbb{Q}}
g^{*}(K_{X}+M_{X})+\Sigma_{X}\sim_{\mathbb{Q}}
f^{*}(K_{Y}+M_{Y})+\Sigma_{Y},
$$
where $M_{W}=g^{-1}(M_{X})$, $\Sigma_{X}$ and $\Sigma_{Y}$ are
exceptional divisors of $g$ and $f$ respectively.

The canonicity of log pairs $(X, M_{X})$ and $(Y, M_{Y})$ implies
the effectiveness of the exceptional divisors $\Sigma_{X}$ and
$\Sigma_{Y}$. However, the effectiveness of $\Sigma_{X}$ and
$\Sigma_{Y}$ implies that
$$
\mathrm{dim}(|km(K_{W}+M_{W})|)=\mathrm{dim}(|g^{*}(km(K_{X}+M_{X}))|)=\mathrm{dim}(|f^{*}(km(K_{Y}+M_{Y}))|)
$$
for $k\gg 0$ if they are not empty and
$$
\phi_{|km(K_{W}+M_{W})|}=\phi_{|g^{*}(km(K_{X}+M_{X}))|}=\phi_{|f^{*}(km(K_{Y}+M_{Y}))|},
$$
which implies the claim.
\end{proof}

By definition, the Kodaira dimension of a movable log pair is a
birational invariant and a non-decreasing function of the
coefficients of the movable boundary.

\begin{definition}
\label{definition:canonical-model} For a given movable log pair
$(X, M_{X})$, a movable log pair $(V, M_{V})$ is called a
canonical model of $(X, M_{X})$ if $M_{V}=\psi(M_{X})$ for a
birational map $\psi:X\dasharrow V$, the divisor $K_{V}+M_{V}$ is
ample, and singularities of $(V, M_{V})$ are canonical.
\end{definition}

The given definition of a canonical model of a movable log pair
coincide with the classical definition of a canonical model in the
case of empty boundary (see \cite{KMM}). The existence of the
canonical model of a movable log pair implies that its Kodaira
dimension equals to the dimension of the variety.

\begin{lemma}
\label{lemma:canonical-model-is-unique} A canonical model of a
movable log pair is unique if it exists.
\end{lemma}

\begin{proof}
Let $(X, M_{X})$ and $(V, M_{V})$ be canonical models such that
$M_{X}=\rho(M_{V})$ for a birational map $\rho:V\dasharrow X$.
Consider a commutative diagram
\[\xymatrix{
&&W\ar@{->}[ld]_{g}\ar@{->}[rd]^{f}&&\\%
&X\ar@{-->}[rr]_{\rho}&&Y&}\] %
such that $W$ is smooth, $g:W\to X$ and $f:W\to Y$ are birational
morphisms. Then
$$
K_{W}+M_{W}\sim_{\mathbb{Q}}
g^{*}(K_{X}+M_{X})+\Sigma_{X}\sim_{\mathbb{Q}}
f^{*}(K_{Y}+M_{Y})+\Sigma_{Y},
$$
where $M_{W}=g^{-1}(M_{X})$, $\Sigma_{X}$ and $\Sigma_{Y}$ are
exceptional divisors of $g$ and $f$ respectively. Then
$$
K_{W}+M_{W}\sim_{\mathbb{Q}}
g^{*}(K_{X}+M_{X})+\Sigma_{X}\sim_{\mathbb{Q}}
f^{*}(K_{V}+M_{V})+\Sigma_{V},
$$
where $M_{W}=g^{-1}(M_{X})=f^{-1}(M_{V})$, and $\Sigma_{X}$ and
$\Sigma_{V}$ are the exceptional divisors of birational morphisms
$g$ and $f$ respectively. The canonicity of the singularities of
the movable log pairs $(X, M_{X})$ and $(V, M_{V})$ implies that
$\Sigma_{X}$ and $\Sigma_{V}$ are effective.

Let $n\in \mathbb{N}$ be a big and divisible enough number such
that $n(K_{W}+M_{W})$, $n(K_{X}+M_{X})$ and $n(K_{V}+M_{V})$ are
Cartier divisors. Then the effectiveness of $\Sigma_{X}$ and
$\Sigma_{V}$ implies
$$
\phi_{|n(K_{W}+M_{W})|}=\phi_{|g^{*}(n(K_{X}+M_{X}))|}=\phi_{|f^{*}(n(K_{V}+M_{V}))|}
$$
and $\rho$ is an isomorphism, because $K_{X}+M_{X}$ and
$K_{V}+M_{V}$ are ample.
\end{proof}

In the case of empty movable boundary the claim of
Lemma~\ref{lemma:canonical-model-is-unique} about the uniqueness
of a canonical modal of an algebraic variety is well known. The
latter implies that all birational automorphisms of a canonical
model are biregular. However, the absence of non-biregular
birational automorphisms is also a property of a birationally
superrigid variety (see
Definition~\ref{definition:superrigidity}). We show later that
Lemma~\ref{lemma:canonical-model-is-unique} explains the
geometrical nature of this phenomenon in the both cases. In the
case of a birationally rigid varieties
Lemma~\ref{lemma:canonical-model-is-unique} is nothing but a
veiled Noether--Fano--Iskovskikh inequality (see \cite{Pu00a}).

\section{Preliminary results.}
\label{sec:3}

Properties of movable log pairs (see
Definition~\ref{definition:movable-log-pair}) reflects birational
geometry of a given variety (see
Lemma~\ref{lemma:canonical-model-is-unique}). Canonical and
terminal singularities are most appropriate classes of
singularities for movable log pairs (see
Remark~\ref{remark:LMMP}). Many geometrical problems can be
translated into the language of movable log pairs. Movable log
pair always can be considered as usual log pairs, and movable
boundaries always can be considered as effective divisors. On the
other hand, we can consider log pairs with both movable and fixed
components (linear systems can have both movable and fixed parts).
Moreover, we can consider log pairs with negative coefficients s
well. We must consider such generalizations due to several
reasons.

For a movable log pair $(X, M_{X})$ and birational morphism
$f:V\to X$, the birationally equivalent log pair $(V, M_{V})$ does
not reflect the properties of the log pair $(X, M_{X})$, but the
log-pullback (see Definition~\ref{definition:log-pull-back}) of
the log pair $(X, M_{X})$ reflects the properties of the log pair
$(X, M_{X})$. However, the log pull back $(V, M^{V})$  of the
movable log pair $(X, M_{X})$ is not necessary a movable log pair
and $M^{V}$ is not necessary an effective divisor. This is the
first reason to consider log pairs with both fixed and movable
components and possibly negative coefficients.

Canonical singularities and centers of canonical singularities
(see Definition~\ref{definition:center}) do not have good
functorial properties when considered apart from the birational
context, but log canonical singularities and centers of log
canonical singularities (see
Definition~\ref{definition:log-center}) have good functorial
properties, and they role in the modern algebraic geometry is very
important (see \cite{KMM}, \cite{Ko91}, \cite{Ko97}, \cite{Sho93},
\cite{Mu01}, \cite{Mu02}, \cite{Pu02a}, \cite{dFEM03}). Log
canonical singu\-la\-ri\-ties and canonical singularities are
related mostly through the log adjunction (see \cite{Co00} and
Theorem~\ref{theorem:log-adjunction}), but also through other ways
(see \cite{Pu02a}). However, log adjunction for movable log pair
can lead to a non-movable log pair. This is another reason to
consider log pairs with both fixed and movable components.

In this section we do not impose any restrictions on boundaries.
In particular, boundaries may not be effective unless otherwise
stated. For simplicity, we assume that log canonical divisors of
all log pairs are $\mathbb{Q}$-Cartier divisors.

\begin{definition}
\label{definition:log-pull-back} A log pair $(V, B^{V})$ is called
a log pull back of a log pair $(X, B_{X})$ with respect to a
birational morphism $f:V\to X$ if we have
$$
B^{V}=f^{-1}(B_{X})-\sum_{i=1}^{n}a(X, B_{X}, E_{i})E_{i}
$$
such that the equivalence
$K_{V}+B^{V}\sim_{\mathbb{Q}}f^{*}(K_{X}+B_{X})$ holds, where
$E_{i}$ is an $f$-exceptional divisor and $a(X, B_{X},
E_{i})\in\mathbb{Q}$.  The rational number $a(X, B_{X}, E_{i})$ is
called a discrepancy of the log pair $(X, B_X)$ in the
$f$-exceptional divisor $E_i$.
\end{definition}

\begin{definition}
\label{definition:log-center} A proper irreducible subvariety
$Y\subset X$ is called a center of log canonical singularities of
the log pair $(X, B_{X})$ if there are a birational morphism
$f:V\to X$ and a not necessary $f$-exceptional divisor $E\subset
V$ such that $E$ is contained in the effective part of the support
of the divisor $\lfloor B^{V}\rfloor$ and $f(E)=Y$.
\end{definition}

\begin{definition}
\label{definition:set-of-log-centers}  The set of centers of log
canonical sin\-gu\-la\-ri\-ties of the log pair $(X, B_{X})$ is
denoted as $\mathbb{LCS}(X, B_{X})$. The set-theoretic union of
all elements in $\mathbb{LCS}(X, B_{X})$ is called a locus of log
canonical singularities of the log pair $(X, B_{X})$, it is
denoted as $LCS(X, B_{X})$.
\end{definition}

\begin{remark}
\label{remark:log-canonical-reduction} Let $H$ be a general
hyperplane section of $X$ and $Z\in \mathbb{LCS}(X, B_{X})$ such
that the inequality $\mathrm{dim}(Z)\ge 1$ holds. Then $Z\cap H\in
\mathbb{LCS}(H, B_{X}\vert_{H})$.
\end{remark}

Let $X$ be a variety and $B_X=\sum_{i=1}^{n}a_{i}B_{i}$ be a
boundary on $X$, where $a_{i}$ is a rational number and $B_{i}$ is
a prime divisor on $X$. Let $f:V\to X$ be a birational morphism
such that $V$ is smooth and the union all $f$-exceptional divisors
and $\cup_{i=1}^{n}f^{-1}(B_{i})$ forms a divisor with simple
normal crossing. The morphism $f$ is called is called a log
resolution of the log pair $(X, B_{X})$. Then the equivalence
$$
K_{Y}+B^{Y}\sim_{\mathbb{Q}} f^{*}(K_{X}+B_{X})
$$
holds, where $(Y, B^{Y})$ is a log pull back of the log pair $(X,
B_{X})$.

\begin{definition}
\label{definition:log-subscheme} Let $\mathcal{I}(X,
B_{X})=f_{*}(\mathcal{O}_{V}(\lceil -B^{V}\rceil))$. Then the
subscheme $\mathcal{L}(X, B_{X})$ associated to the ideal sheaf
$\mathcal{I}(X, B_{X})$ is called a log canonical singularity
subscheme of the log pair $(X, B_{X})$.
\end{definition}

Note, that by definition we have $\mathrm{Supp}(\mathcal{L}(X,
B_{X}))=LCS(X, B_{X})\subset X$. The following result is the
Shokurov vanishing theorem (see \cite{Sho93}, \cite{Am99}).

\begin{theorem}
\label{theorem:vanishing-of-Shokurov} Let $(X, B_{X})$ be a log
pair,  and $H$ be a nef and big divisor on $X$ such that the
boundary $B_{X}$ is effective, and $D=K_{X}+B_{X}+H$ is a Cartier
divisor. Then the cohomology group $H^{i}(X, \mathcal{I}(X,
B_{X})\otimes D)$ vanishes for $i>0$.
\end{theorem}

\begin{proof}
Let $f:W\longrightarrow X$ be a log resolution of $(X,B_X)$. Then
$$
R^{i}f_{*}(f^{*}(K_{X}+B_{X}+H)+\lceil -B^{W}\rceil)=0
$$
for $i>0$ by the Kawamata-Viehweg vanishing (see \cite{Ka82},
\cite{Vi82}, \cite{KMM}). The degeneration of the
local--to--global spectral sequence and
$$R^{0}f_{*}(f^{*}(K_{X}+B_{X}+H)+\lceil -B^{W}\rceil)=\mathcal{I}(X, B_{X})\otimes D$$
imply that for all $i\ge 0$ we have
$$H^{i}(X, \mathcal{I}(X, B_{X})\otimes D)=H^{i}(W, f^{*}(K_{X}+B_{X}+H)+\lceil -B^{W}\rceil),$$
but $H^{i}(W, f^{*}(K_{X}+B_{X}+H)+\lceil -B^{W}\rceil)=0$ for
$i>0$ by the Kawamata-Viehweg vanishing.
\end{proof}

Consider the following two application of
Theorem~\ref{theorem:vanishing-of-Shokurov}, which are special
cases of a more general result in \cite{Ch01b} (see \cite{ChPa02},
\cite{EM01}, \cite{dFEM03}).

\begin{lemma}
\label{lemma:quadric-surface} Let $V=\mathbb{P}^{1}\times
\mathbb{P}^{1}$ and $B_{V}$ be an effective boundary on $V$ of
bi-degree $(a, b)$ such that $a$ and $b\in\mathbb{Q}\cap [0, 1)$.
Then $\mathbb{LCS}(V, B_{V})=\emptyset$.
\end{lemma}

\begin{proof}
Let $B_{V}=\sum_{i=1}^{k}a_{i}B_{i}$, where $a_{i}$ is a positive
rational number, and $B_{i}$ is an irreducible reduced curve on
the surface $V$. Intersecting the boundary $B_{V}$ with the
rulings of $V$ we get the inequality $a_{i}<1$. Thus the set
$\mathbb{LCS}(V, B_{V})$ does not contains curves on $V$.

Suppose that the set $\mathbb{LCS}(V, B_{V})$ contains a point
$O$. Take a divisor $H\in\mathrm{Pic}(V)\otimes\mathbb{Q}$  of
bi-degree $(1-a, 1-b)$. Then the divisor $H$ is ample. Moreover,
there is a divisor
$$D\sim_{\mathbb{Q}} K_{V}+B_{V}+H$$
such that $D$ is a Cartier divisor and
$H^{0}(\mathcal{O}_{V}(D))=0$. On the other hand, the map
$$
H^{0}(\mathcal{O}_{V}(D))\to H^{0}(\mathcal{O}_{\mathcal{L}(V,
B_{V})}(D))
$$
is surjective by Theorem~\ref{theorem:vanishing-of-Shokurov},
which is a contradiction.
\end{proof}

\begin{lemma}
\label{lemma:cubics} Let $V\subset\mathbb{P}^{n}$ be a smooth
hypersurface of degree $k<n$, and $B_{V}$ be an effective boundary
on $V$ such that $B_{V}\equiv rH$, where $r\in\mathbb{Q}\cap [0,
1)$, and $H$ is a hyperplane section of the hypersurface
$V\subset\mathbb{P}^{n}$. Then $\mathbb{LCS}(V, B_{V})=\emptyset$.
\end{lemma}

\begin{proof}
Suppose that the set $\mathbb{LCS}(V, B_{V})$ contains a
subvariety $Z\subset V$. Then $\mathrm{dim}(Z)=0$ by Theorem~2 in
\cite{Pu95} (see Lemma~3.18 in \cite{Ch03b}). Therefore the set
$\mathbb{LCS}(V, B_{V})$ contains only closed points of the
hypersurface $V$. In particular, the support of the scheme
$\mathcal{L}(V, B_{V})$ is zero-dimensional and
$H^{0}(\mathcal{O}_{\mathcal{L}(V,\, B_{V})})\ne 0$.

Note, that $K_{V}+B_{V}+(1-r)H\equiv (k-n)H$ and
$H^{0}(\mathcal{O}_{V}((k-n)H))=0$, because the inequality $k<n$
holds. However, Theorem~\ref{theorem:vanishing-of-Shokurov}
implies the surjectivity
$$
H^{0}(\mathcal{O}_{V}((k-n)H))\to
H^{0}(\mathcal{O}_{\mathcal{L}(V, B_{V})}((k-n)H))\to 0,
$$
which is a contradiction, because
$H^{0}(\mathcal{O}_{\mathcal{L}(V,\,
B_{V})}((k-n)H))=H^{0}(\mathcal{O}_{\mathcal{L}(V,\, B_{V})})$.
\end{proof}

\begin{example}
\label{example:lemma-cubics-fails-when-k-equals-n} Let
$V\subset\mathbb{P}^{n}$ be a smooth hypersurface
$$
x_{0}^{k}=\sum_{i=1}^{n}x_{i}^{k}\subset\mathbb{P}^{n}\cong\mathrm{Proj}(\mathbb{C}[x_{0},\ldots,x_{n}]),
$$
and $B_{V}={\frac{n-1}{k}}H$, where $H$ is a hyperplane section of
the hypersurface $V$ that is cut by the equation $x_{0}=x_{1}$.
Then the hypersurface $V$ is smooth and the set $\mathbb{LCS}(V,
B_{V})$ consist of a single point $(1:1:0:\ldots:0)\in
V\subset\mathbb{P}^{n}$.
\end{example}

The arguments of the proofs of Lemmas~\ref{lemma:quadric-surface}
and \ref{lemma:cubics} can be applied in much more general
situation. Namely, for a given Cartier divisor $D$ on the variety
$X$, let us consider the exact sequence of sheaves
$$
0\to\mathcal{I}(X, B_{X})\otimes D\to\mathcal{O}_{X}(D)\to\mathcal{O}_{\mathcal{L}(X,\, B_{X})}(D)\to 0,%
$$
and the corresponding exact sequence of cohomology groups. Now
Theorem~\ref{theorem:vanishing-of-Shokurov} implies the following
two connectedness results (see \cite{Sho93}).

\begin{theorem}
\label{theorem:global-connectedness} Let $(X, B_{X})$ be a log
pair, and let $B_{X}$ be an effective boundary such that the
divisor $-(K_{X}+B_{X})$ is nef and big. Then the locus $LCS(X,
B_{X})$ is connected.
\end{theorem}

\begin{theorem}
\label{theorem:local-connectedness} Let $(X, B_{X})$ be a log
pair, $B_{X}$ be an effective boundary, $g:X\longrightarrow Z$ be
morphism with connected fibers such that $-(K_{X}+B_{X})$ is
$g$-nef and $g$-big. Then $LCS(X, B_{X})$ is con\-nected in a
neighborhood of each fiber of $g$.
\end{theorem}

Similarly, one can prove the following  result, which is
Theorem~17.4 in \cite{Ko91}.

\begin{theorem}
\label{theorem:general-connectedness} Let $g:X\to Z$ be a morphism
with connected fibers, $D=\sum_{i\in\, I} d_{i}D_{i}$ be a divisor
on $X$, $h:V\to X$ be a resolution of singularities of the variety
$X$  such that the union of all $h$-exceptional divisors and
$\cup_{i\in I}h^{-1}(D_{i})$ is a simple normal crossing divisor,
the divisor $-(K_{X}+D)$ is $g$-nef and $g$-big, and the
inequality $\mathrm{codim}(g(D_{i})\subset Z)\geq 2$ holds
whenever $d_{i}<0$. For any divisor $E\subset V$ let
$a(E)\in\mathbb{Q}$ such that the equivalence
$$
K_{V}\sim_\mathbb{Q} f^{*}(K_{X}+D)+\sum_{E\subset V} a(E)E
$$
holds. Then $\cup_{a(E)\le -1}E$ is connected in the neighborhood
of every fiber of $g\circ h$.
\end{theorem}

\begin{proof}
Let $f=g\circ h$, $A=\sum_{a(E)>-1}E$, and $B=\sum_{a(E)\le -1}E$.
Then
$$
\lceil A\rceil-\lfloor B\rfloor\sim_{\mathbb{Q}}
K_{V}-h^{*}(K_{X}+D)+\lbrace-A\rbrace+\lbrace B\rbrace
$$
and $R^{1}f_{*}\mathcal{O}_{V}(\lceil A\rceil-\lfloor B\rfloor)=0$
by the Kawamata--Viehweg vanishing. Hence the map
$$
f_{*}\mathcal{O}_{V}(\lceil A\rceil)\to f_{*}\mathcal{O}_{\lfloor
B\rfloor}(\lceil A\rceil)
$$
is surjective. Every component of $\lceil A\rceil$ is either
$h$-ex\-cep\-ti\-onal or a proper transform of a divisor $D_{j}$
with $d_{j}<0$. Thus $h_{*}(\lceil A\rceil)$ is $g$-exceptional
and $f_{*}\mathcal{O}_{V}(\lceil A\rceil)=\mathcal{O}_{Z}$. So the
map
$$
\mathcal{O}_{Z}\to f_{*}\mathcal{O}_{\lfloor B\rfloor}(\lceil A\rceil)%
$$
is surjective, which implies the connectedness of $\lfloor
B\rfloor$ in a neighborhood of every fiber of morphism $f$,
because $\lceil A\rceil$ is effective and has no common component
with $\lfloor B\rfloor$.
\end{proof}

We defined the notion of a center of canonical singularities in
Definition~\ref{definition:center} for a movable log pair.
However, we did not use the movability of a boundary in
Definitions~\ref{definition:center}, and we can consider centers
of canonical singularities of any log pair.

\begin{theorem}
\label{theorem:log-adjunction} Let $(X, B_{X})$ be a log pair, $Z$
be an element in $\mathbb{CS}(X, B_{X})$, $H$ be an effective
irreducible Cartier divisor on $X$ such that $Z\subset H$, $X$ and
$H$ are smooth at the generic point of $Z$, $H$ is not a component
of $B_{X}$, and $B_{X}$ is effective. Then $\mathbb{LCS}(H,
B_{X}\vert_{H})\ne\emptyset$.
\end{theorem}

\begin{proof}
Let $f:W\to X$ be a log resolution of $(X, B_{X}+H)$. Put ${\hat
H}=f^{-1}(H)$. Then
$$K_{W}+{\hat H}\sim_{\mathbb{Q}} f^{*}(K_{X}+B_{X}+H)+\sum_{E\ne {\hat H}} a(X, B_{X}+H, E)E$$
and by assumption we have $\{ Z, H\}\subset LCS(X, B_{X}+H)$.
Therefore applying Theorem~\ref{theorem:general-connectedness} to
the log pullback of $(X, B_{X}+H)$ on $W$, we get ${\hat H}\cap
E\ne\emptyset$ for some divisor $E\ne {\hat H}$ on the variety $W$
such that $f(E)=Z$ and $a(X, B_{X}, E)\le -1$. Now the
equivalences
$$K_{\hat H}\sim (K_{W}+{\hat H})\vert_{\hat H}\sim_{\mathbb{Q}} f\vert_{\hat H}^{*}(K_{H}+B_{X}\vert_{H})+\sum_{E\ne {\hat H}} a(X, B_{X}+H, E)E\vert_{\hat H}$$
imply the claim.
\end{proof}

\begin{corollary}
\label{corollary:log-adjunction-smooth-points} Let $(X, M_{X})$ be
a movable log pair, $O$ be a smooth point of $X$, $H_{i}$ be a
general hyperplane section of $X$ passing through the point $O$
for $i=1,\ldots,k\le\mathrm{dim}(X)-2$ such that
$O\in\mathbb{CS}(X, M_{X})$, $M_{X}$ is effective, and
$\mathrm{dim}(X)\ge 3$. Then $O\in\mathbb{LCS}(S, M_{S})$, where
$S=\cap_{i=1}^{k}H_{i}$ and $M_{S}=M_{X}\vert_{S}$.
\end{corollary}

It should be pointed out that Theorem~\ref{theorem:log-adjunction}
is a special case of a general phenomenon, which is known as log
adjunction (see \cite{Ko91}, \cite{Co00}). In particular, simple
modification of the proof of Theorem~\ref{theorem:log-adjunction}
implies the following result.

\begin{corollary}
\label{corollary:log-adjunction-double-points} Let $(X, M_{X})$ be
a movable log pair, $O$ be an isolated hypersurface singular point
of the variety $X$, and $H_{i}$ be a general hyperplane section of
$X$ passing through the point $O$ for
$i=1,\ldots,k\le\mathrm{dim}(X)-2$ such that $O\in\mathbb{CS}(X,
M_{X})$, the boundary $M_{X}$ is effective, and
$\mathrm{dim}(X)\ge 3$. Then $O\in\mathbb{LCS}(S, M_{S})$, where
$S=\cap_{i=1}^{k}H_{i}$ and $M_{S}=M_{X}\vert_{S}$.
\end{corollary}

The following result is a Theorem~3.1 in \cite{Co00}, which gives
the shortest proof of the main result of \cite{IsMa71} modulo
Theorem~\ref{theorem:log-adjunction} (see \cite{Co00}).

\begin{theorem}
\label{theorem:Corti} Suppose that $\mathrm{dim}(X)=2$, the
boundary $B_{X}$ is effective and movable, and there is a smooth
point $O\in X$ such that $O\in\mathbb{LCS}(X,
(1-a_{1})\Delta_{1}+(1-a_{2})\Delta_{2}+M_{X})$, where
$\Delta_{1}$ and $\Delta_{2}$ are smooth curves on $X$
intersecting normally at $O$, and $a_{1}$ and $a_{2}$ are
arbitrary non-negative rational numbers. Then we have
$$
\mathrm{mult}_{O}(B_{X}^{2})\ge\left\{\aligned
&4a_{1}a_{2}\ {\text {if}}\ a_{1}\le 1\ {\text {or}}\ a_{2}\le 1\\
&4(a_{1}+a_{2}-1)\ {\text {if}}\ a_{1}>1\ {\text {and}}\ a_{2}>1.\\
\endaligned
\right.
$$
\end{theorem}

Most applications of Theorem~\ref{theorem:Corti} use the
simplified version Theorem~\ref{theorem:Corti} (see
Corollary~\ref{corollary:Corti-simple-form}) that involves only
movable boundary. Moreover, Theorem~\ref{theorem:Corti} was
created in order to be applied to movable log pairs. However, the
proof of Theorem~\ref{theorem:Corti} in \cite{Co00} is inductive
by the number of blow ups required to obtain the appropriate
negative discrepancy. It is easy to see that the inductive proof
of Theorem~\ref{theorem:Corti} is much easier to apply when we
have nonmovable components of the boundary. In certain sense the
main difficulty in the proof of Theorem~\ref{theorem:Corti} is to
find the right form of Theorem~\ref{theorem:Corti}, which is
suitable for the inductive proof. On the other hand,
Theorem~\ref{theorem:Corti} with nontrivial nonmovable components
of the boundary has nice higher-dimensional applications (see
\cite{Ch00b}, \cite{Ch03b}). More general approach to
Theorem~\ref{theorem:Corti} was found in \cite{dFEM03}, where an
analog of Theorem~\ref{theorem:Corti} was used to prove the
generalization of the main inequality of \cite{Pu02a}. Note, that
Theorem~2.1 in \cite{dFEM03} is a generalization of
Theorem~\ref{theorem:Corti} in the case when the nonmovable part
of the boundary consists of a single component. However, such
weaken version of Theorem~\ref{theorem:Corti} is not suitable for
some applications (see \cite{Ch00b}).

The following result is a special case of Theorem~0.1 in
\cite{dFEM02}.

\begin{corollary}
\label{corollary:Corti-simple-form} Let $H$ be a surface, $O$ be a
smooth point on $H$, and $M_{H}$ be an effective movable boundary
on $H$ such that $O\in\mathbb{LCS}(H, M_{H})$. Then the inequality
$\mathrm{mult}_{O}(M_{H}^{2})\geq 4$ holds and the equality
$\mathrm{mult}_{O}(M_{H}^{2})=4$ implies
$\mathrm{mult}_O(M_{H})=2$.
\end{corollary}

The following result is due to \cite{Pu00a}.

\begin{theorem}
\label{theorem:Iskovskikh} Let $X$ be a variety, $M_{X}$ be an
effective movable boundary on $X$, and $O$ be a smooth point of
$X$ such that $O\in\mathbb{CS}(X, M_{X})$ and $\mathrm{dim}(X)\ge
3$. Then $\mathrm{mult}_{O}(M_{X}^{2})\ge 4$ and the equality
$\mathrm{mult}_{O}(M_{X}^{2})=4$ implies
$\mathrm{mult}_{O}(M_{X})=2$ and $\mathrm{dim}(X)=3$.
\end{theorem}

\begin{proof}
The claim is implied by
Corollaries~\ref{corollary:log-adjunction-smooth-points} and
\ref{corollary:Corti-simple-form}.
\end{proof}

The proof of Theorem~\ref{theorem:Iskovskikh} in \cite{Pu00a} is
elementary but technical, which is valid even over fields of
positive characteristic. The proof in \cite{Pu00a} and the proof
in \cite{Co00} does not explain the geometrical nature of
Theorem~\ref{theorem:Iskovskikh}, which is pointed out in
\cite{Co95} and requires the following well known result (see
\cite{Ko91}).

\begin{lemma}
\label{lemma:terminal-extraction} Let $X$ be a smooth 3-fold, $O$
be a point on $X$, and $M_{X}$ be an effective movable boundary on
the variety $X$ such that the singularities of the log pair $(X,
M_{X})$ are canonical, and $O\in\mathbb{CS}(X, M_{X})$. Then there
is a birational morphism $f:V\to X$ such that the 3-fold $V$ has
$\mathbb{Q}$-factorial terminal singularities, the morphism $f$
contracts exactly one divisor $E$, $f(E)=O$, and
$K_{V}+M_{V}\sim_{\mathbb{Q}} f^{*}(K_{X}+M_{X})$, where
$M_{V}=f^{-1}(M_{X})$.
\end{lemma}

\begin{proof}
There are finitely many divisorial discrete valuations $\nu$ of
the field of rational functions of $X$ whose center on $X$ is the
point $O$ and whose discrepancy $a(X,M_X,\nu)$ is non-positive,
because $(X, M_X)$ has canonical singularities. Therefore we may
consider a birational morphism $g:W\to X$ such that $W$ is smooth,
$g$ contracts $k$ divisors,
$$
K_{W}+M_{W}\sim_\mathbb{Q}
g^{*}(K_{X}+M_{X})+\sum_{i=1}^{k}a_{i}E_{i},
$$
movable log pair $(W, M_{W})$ has canonical singularities, and the
set $\mathbb{CS}(W, M_{W})$ does not contain subvarieties of
$\cup_{i=1}^{k}E_{i}$, where $M_{W}=g^{-1}(M_{X})$, $g(E_{i})=O$,
and $a_{i}\in \mathbb{Q}$. Applying the relative Log Minimal Model
Program (see \cite{KMM}) to the log pair $(W, M_{W})$ over $X$, we
may assume that the 3-fold $W$ has terminal $\mathbb{Q}$-factorial
singularities and
$$
K_{W}+M_{W}\sim_\mathbb{Q} g^{*}(K_{X}+M_{X})
$$
because the singularities of the movable log pair $(X, M_{X})$ are
canonical. Now applying the relative Minimal Model Program to the
variety $W$ over $X$, we get the necessary 3-fold and birational
morphism.
\end{proof}

\begin{remark}
\label{remark:terminal-extraction} Let $X$ be a smooth variety,
$O$ be a point of $X$, and $f:V\to X$ be a birational morphism
such that $V$ has terminal $\mathbb{Q}$-factorial singularities,
$f$ contracts a single exceptional divisor $E$, and $f(E)=O$. Then
there is a movable log pair $(X, M_{X})$ such that the boundary
$M_{X}$ is effective, the singularities of the log pair $(X,
M_{X})$ are canonical, $K_{V}+M_{V}\sim_{\mathbb{Q}}
f^{*}(K_{X}+M_{X})$, and $O\in\mathbb{CS}(X, M_{X})$, where
$M_{V}=f^{-1}(M_{X})$.
\end{remark}

The following result is conjectured in \cite{Co95} and proved in
\cite{Kaw01}.

\begin{theorem}
\label{theorem:Kawakita} Let $X$ be a smooth 3-fold, $O$ be a
point of $X$, and $f:V\to X$ be a birational morphism such that
the singularities of $V$ are terminal and $\mathbb{Q}$-factorial,
$f$ contracts a single divisor $E\subset V$, and $f(E)=O$. Then
$f$ is a weighted blow up at the point $O$ with weights $(1,K,N)$
in suitable local coordinates on $X$, where $K$ and $N$ are
coprime naturals.
\end{theorem}

Actually, Theorem~\ref{theorem:Iskovskikh} was proved in
\cite{Co95} modulo Theorem~\ref{theorem:Kawakita} in the following
way, which explains the geometrical nature of
Theorem~\ref{theorem:Iskovskikh}.

\begin{proposition}
\label{proposition:weighted-intersections} Let $X$ be a smooth
3-fold, $O$ be a point of $X$, and $M_{X}$ be an effective movable
boundary on $X$, and $f:V\to X$ be a weighted blow up of $O$ with
weights $(1,K,N)$ in suitable local coordinates on $X$ such that
$$
K_{V}+M_{V}\sim_{\mathbb{Q}} f^{*}(K_{X}+M_{X}),
$$
where natural numbers $K$ and $N$ are coprime and
$M_{V}=f^{-1}(M_{X})$. Then
$$
\mathrm{mult}_{O}(M_{X}^{2})\geq {\frac {(K+N)^{2}} {KN}}=4+{\frac
{(K-N)^{2}} {KN}}\geq 4,
$$
where $K=N$ implies that $f$ is a standard blow up of $O$ and
$\mathrm{mult}_{O}(M_{X})=2$.
\end{proposition}

\begin{proof}
Let $E\subset V$ be an $f$-ex\-cep\-ti\-o\-nal divisor. Then
$$
K_{V}\sim_{\mathbb{Q}} f^{*}(K_{X})+(N+K)E
$$
and $M_{V}\sim_{\mathbb{Q}} f^{*}(M_{X})+mE$ for some $m\in
{\mathbb{Q}}_{>0}$. Thus $m=K+N$. Now intersecting the effective
cycle $M_{X}^{2}$ with a general hyperplane section of $X$ passing
through $O$, we obtain the inequality
$\mathrm{mult}_{O}(M_{X}^{2})\ge m^{2}E^{3}={\frac {(K+N)^{2}}
{KN}}$.
\end{proof}

The following application of Theorem~\ref{theorem:log-adjunction}
is Theorem 3.10 in \cite{Co00}.

\begin{theorem}
\label{theorem:double-point} Let $X$ be a variety, $O$ be an
ordinary double point of $X$, and $B_{X}$ an effective boundary on
the variety $X$ such that $O\in\mathbb{CS}(X, B_{X})$, $B_{X}$ is
a $\mathbb{Q}$-Cartier divisor, and $\mathrm{dim}(X)\ge 3$. Then
$\mathrm{mult}_{O}(B_{X})\geq 1$, and $\mathrm{mult}_{O}(B_{X})=1$
implies $\mathrm{dim}(X)=3$, where the positive rational number
$\mathrm{mult}_{O}(B_{X})$ is defined through the standard blow up
of $O$.
\end{theorem}

\begin{proof}
We may assume that $X$ is a 3-fold due to
Corollary~\ref{corollary:log-adjunction-double-points}. Let
$f:W\to X$ be a blow up at the point $O$ and $E$ be an
$f$-exceptional divisor. Then
$$
K_{W}+B_{W}\sim_\mathbb{Q}
f^{*}(K_{X}+B_{X})+(1-\mathrm{mult}_{O}(B_{X}))E,
$$
where $B_{W}=f^{-1}(B_{X})$. Suppose that the inequality
$\mathrm{mult}_{O}(B_{X})<1$ holds. Then there is a proper
subvariety $Z\subset Q$ such that $Z\in\mathbb{CS}(W, B_{W})$.
Hence
$$
\mathbb{LCS}(E, B_{W}\vert_{E})\ne\emptyset
$$
by Theorem~\ref{theorem:log-adjunction}, which is impossible by
Lemma~\ref{lemma:quadric-surface}, because
$E\cong\mathbb{P}^{1}\times\mathbb{P}^{1}$.
\end{proof}

\begin{proposition}
\label{proposition:non-simple-double-point} Let $X$ be a variety,
$B_{X}$ be an effective boundary on $X$, and $O$ be an isolated
singular point on $X$ such that $X$ is locally given by the
equation $y^{3}=\sum_{i=1}^{\mathrm{dim}(X)}x_{i}^{2}$ in the
neighborhood of $O$, the boundary $B_{X}$ is a
$\mathbb{Q}$-Cartier divisor on $X$, $O\in\mathbb{CS}(X, B_{X})$,
and $\mathrm{dim}(X)\ge 4$. Then $\mathrm{mult}_{O}(B_{X})>1$,
where $\mathrm{mult}_{O}(B_{X})\in\mathbb{Q}$ is defined naturally
by means of the standard blow up of the point $O$.
\end{proposition}

\begin{proof}
The claim is implied by
Corollary~\ref{corollary:log-adjunction-double-points} and
Theorem~\ref{theorem:double-point}.
\end{proof}

\begin{theorem}
\label{theorem:tripple-point} Let $X$ be a variety of dimension
$n\ge 4$, $B_{X}$ be an effective boundary on the variety $X$, $O$
be an ordinary triple point\,\footnote{Namely, the point $O$ is an
isolated hypersurface singular point of $X$ such that the
projectivization of the tangent cone to $X$ at the point $O$ is a
smooth hypersurface in $\mathbb{P}^{n-1}$ of degree $3$} of the
variety $X$ such that $O\in\mathbb{CS}(X, B_{X})$, and the
boundary $B_{X}$ is a $\mathbb{Q}$-Cartier divisor on $X$. Then
the inequality $\mathrm{mult}_{O}(B_{X})\geq 1$ holds, and
$\mathrm{mult}_{O}(B_{X})=1$ implies $n=4$, where the rational
number $\mathrm{mult}_{O}(B_{X})$ is defined naturally through the
standard blow up of the point $O$.
\end{theorem}

\begin{proof}
Let $f:W\to X$ be a blow up of the point $O$. Then
$$
K_{W}+B_{W}\sim_\mathbb{Q}
f^{*}(K_{X}+B_{X})+(n-3-\mathrm{mult}_{O}(B_{X}))E,
$$
where $B_{W}=f^{-1}(B_{X})$ and $E=f^{-1}(O)$. Suppose that
$\mathrm{mult}_{O}(B_{X})<n-3$. Then there is a subvariety
$Z\subset E$ such that
$$
Z\in\mathbb{CS}(W,
B_{W}-(n-3-\mathrm{mult}_{O}(B_{X}))E)\subseteq\mathbb{CS}(W,
B_{W}),
$$
and the inequalities $n>4$ and $\mathrm{mult}_{O}(B_{X})\le 1$
imply that
$$
\mathbb{CS}(W,
B_{W}-(n-3-\mathrm{mult}_{O}(B_{X}))E)\subseteq\mathbb{CS}(W,
\lambda B_{W})
$$
for some positive rational number $\lambda<1$. In particular,
$\mathbb{LCS}(E, B_{W}\vert_{E})\ne\emptyset$ in the case when
$\mathrm{mult}_{O}(B_{X})<1$ by
Theorem~\ref{theorem:log-adjunction}. Moreover, we have
$\mathbb{LCS}(E, \lambda B_{W}\vert_{E})\ne\emptyset$ in the case
when $\mathrm{mult}_{O}(B_{X})\le 1$ and $n>4$. Therefore, in both
cases we proved the claim that contradicts to
Lemma~\ref{lemma:cubics}.
\end{proof}

It is easy to see that Theorems~\ref{theorem:double-point} and
\ref{theorem:tripple-point} are special cases of the following
general result, which is left without a proof, because its proof
is very similar to the proof of
Theorem~\ref{theorem:tripple-point}.

\begin{theorem}
\label{theorem:mulpitle-point} Let $X$ be a variety of dimension
$n$, $B_{X}$ be an effective boundary on the variety $X$, and $O$
be an ordinary singular point\footnote{Namely, the point $O$ is an
isolated hypersurface singular point on $X$ such that the
projectivization of the tangent cone to $X$ at the point $O$ is a
smooth hypersurface in $\mathbb{P}^{n-1}$.} of multiplicity $k$
such that $O\in\mathbb{CS}(X, B_{X})$, the inequality $n>k$ holds,
and $B_{X}$ is a $\mathbb{Q}$-Cartier divisor. Then
$\mathrm{mult}_{O}(B_{X})\geq 1$, and the equality
$\mathrm{mult}_{O}(B_{X})=1$ implies $n=k+1$, where
$\mathrm{mult}_{O}(B_{X})\in\mathbb{Q}$ is defined naturally
through the standard blow up of the point $O$.
\end{theorem}

\begin{corollary}
\label{corollary:mulpitle-point} Let $f:V\to X$ be a birational
morphism, $O$ be an ordinary singular point of the variety $X$ of
multiplicity $\mathrm{dim}(X)-1$ such that $X$ and $V$ have
terminal $\mathbb{Q}$-factorial singularities, the morphism $f$
contracts a single divisor $E$, and $f(E)=O$. Then $f$ is a
standard blow up of the point $O$.
\end{corollary}

\section{The Noether--Fano--Iskovskikh inequality.}
\label{sec:4}

In this chapter we consider the Noether--Fano--Iskovskikh
inequality and give two generalization of this inequality. Let $X$
be a Fano variety with terminal $\mathbb{Q}$-factorial
singularities such that $\mathrm{Pic}(X)\cong\mathbb{Z}$. For
example, we can always substitute $X$ by the variety that
satisfies all conditions of Theorems~\ref{theorem:main} or
\ref{theorem:third} (see Lemma~\ref{lemma:factoriality} and
Remark~\ref{remark:factoriality}). We assume that all movable
boundaries are effective. The following result is due to
\cite{Co95}, but its special cases can be found in \cite{Ne71},
\cite{Fa15}, \cite{Fa47}, \cite{Ma66}, \cite{Ma67}, \cite{IsMa71},
\cite{Is80b}, \cite{Sa80} and \cite{Sa82}.

\begin{theorem}
\label{theorem:Nother-Fano-inequality} Suppose that every movable
log pair $(X, M_{X})$ such that $K_{X}+M_{X}\sim_\mathbb{Q} 0$ has
canonical singularities. Then the Fano variety $X$ is birationally
superrigid.
\end{theorem}

\begin{proof}
Suppose that $X$ is not birationally superrigid. Let
$\rho:X\dasharrow Y$ be a birational map such that the rational
map $\rho$ is not biregular and either $Y$ is a Fano variety of
Picard rank $1$ with terminal $\mathbb{Q}$-factorial singularities
or there is a fibration $\tau:Y\to Z$ whose generic fiber has
Kodaira dimension $-\infty$. We may assume that

Suppose that we have a fibration $\tau:Y\to Z$ such that the
generic fiber of $\tau$ is a variety of Kodaira dimension
$-\infty$. Take a very ample divisor $H$ on $Z$ and some positive
rational number $\mu$. Put $M_{Y}=\mu|\tau^{*}(H)|$ and
$M_{X}=\mu\rho^{-1}(|\tau^{*}(H)|)$. Then we have
$$
\kappa(X, M_{X})=\kappa(Y, M_{Y})=-\infty%
$$
by construction. Choose $\mu$ such that $M_{X}\sim_{\mathbb{Q}}
-K_{X}$. Then the singularities of $(X, M_{X})$ are not canonical,
because otherwise $\kappa(X, M_{X})=0$. Therefore we get a
contradiction with our initial assumption.

Suppose that $Y$ is a $\mathbb{Q}$-factorial terminal Fano variety
of Picard rank $1$. Take a positive rational number $\mu$. Let
$M_{Y}={\frac{\mu}{n}}|-nK_{Y}|$ and $M_{X}=\rho^{-1}(M_{Y})$ for
$n\gg 0$. Then
$$
\kappa(X, M_{X})=\kappa(Y, M_{Y})=\left\{\aligned
&\mathrm{dim}(Y)\ \text{for}\ \mu>1\\
&0\ \text{for}\ \mu=1\\
&-\infty\ \text{for}\ \mu<1\\
\endaligned
\right.
$$
by construction. Choose $\mu$ such that $M_{X}\sim_{\mathbb{Q}}
-K_{X}$. Then the singularities of the movable log pair $(X,
M_{X})$ are canonical by assumption. Hence $\kappa(X, M_{X})=0$
and $\mu=1$.

Let us consider a commutative diagram
\[\xymatrix{
&&W\ar@{->}[ld]_{g}\ar@{->}[rd]^{f}&&\\%
&X\ar@{-->}[rr]_{\rho}&&Y&}\] %
such that $W$ is smooth, $g:W\to X$ and $f:W\to Y$ are birational
morphisms. Then
$$
\sum_{j=1}^{k}a(X, M_{X}, F_{j})F_{j}\sim_{\mathbb{Q}}\sum_{i=1}^{l}a(Y, M_{Y}, G_{i})G_{i},%
$$
where $G_{i}$ is an $g$-exceptional divisor and $F_{j}$ is an
$f$-exceptional divisor. We may assume that $k=l$, because $X$ and
$Y$ are $\mathbb{Q}$-factorial and have Picard rank $1$. Every
$a(X, M_{X}, F_{j})$ is non-negative and every $a(Y, M_{Y},
G_{i})$ is positive, because $(Y, M_{Y})$ has terminal
singularities by assumption. Thus it follows from Lemma~2.19 in
\cite{Ko91} that
$$\sum_{j=1}^{k}a(X, M_{X}, F_{j})F_{j}=\sum_{i=1}^{k}a(Y, M_{Y}, G_{i})G_{i},%
$$
which implies that the singularities of the log pair $(X, M_{X})$
are terminal.
\par
Now take $\mu>1$ such that the singularities of the log pairs $(X,
M_{X})$ and $(Y, M_{Y})$ are still terminal (see
Remark~\ref{remark:flexability-terminality}). Then both log pairs
$(X, M_{X})$ and $(Y, M_{Y})$ are canonical models. Thus $\rho$ is
an isomorphism by Lemma~\ref{lemma:canonical-model-is-unique},
which contradicts to our assumption and concludes the proof.
\end{proof}

The roots of Theorem~\ref{theorem:Nother-Fano-inequality} can be
found in \cite{Ne71}, \cite{Fa15} and \cite{Fa47}. In
two-dimensional case an analog of
Theorem~\ref{theorem:Nother-Fano-inequality} is proved in
\cite{Ma66}, \cite{Ma67}, in the three-dimensional case an analog
of Theorem~\ref{theorem:Nother-Fano-inequality} is proved in
\cite{IsMa71}, \cite{Is80b}.

\begin{corollary}
\label{corollary:Nother-Fano} Suppose that $X$ is not birationally
superrigid. Then there is a movable log pair $(X, M_{X})$ such
that the divisor $-(K_{X}+M_{X})$ is ample and $\mathbb{CS}(X,
M_{X})\ne\emptyset$.
\end{corollary}

The following two generalizations of
Theorem~\ref{theorem:Nother-Fano-inequality} are due to
\cite{Ch00a}.

\begin{theorem}
\label{theorem:Nother-Fano-inequality-elliptic-case} Let
$\rho:V\dasharrow X$ be birational map such that there is a
morphism $\tau:V\to Z$ whose generic fiber is a smooth elliptic
curve, and let $D$ be a very ample divisor on the variety $Z$ and
$\mathcal{D}=|\tau^{*}(D)|$. Put $\mathcal{M}=\rho(\mathcal{D})$
and $M_{X}=\gamma\mathcal{M}$, where $\gamma$ is a positive
rational number such that $K_{X}+\gamma M_{X}\sim_\mathbb{Q} 0$.
Then $\mathbb{CS}(X, M_{X})\ne\emptyset$.
\end{theorem}

\begin{proof}
Suppose that the set $\mathbb{CS}(X, M_{X})$ is empty. Then the
singularities of the movable log pair $(X, M_{X})$ are terminal.
In particular, for some rational number $\epsilon>\gamma$ the
movable log pair $(X, \epsilon\mathcal{M})$ is a canonical model
(see Remark~\ref{remark:flexability-terminality}). In particular,
the equality
$$
\kappa(X, \epsilon\mathcal{M})=\mathrm{dim}(X)
$$
holds. On the other hand, the log pairs $(X, \epsilon\mathcal{M})$
and $(V, \epsilon\mathcal{D})$ are  birationally equivalent and
have the same Kodaira dimensions. However, by construction
$$
\kappa(V, \epsilon\mathcal{D})\le
\mathrm{dim}(Z)=\mathrm{dim}(X)-1,
$$
which is a contradiction.
\end{proof}

\begin{theorem}
\label{theorem:Nother-Fano-inequality-Fano-case} Let
$\rho:V\dasharrow X$ be a birational map such that $V$ is a Fano
variety with canonical singularities. Put $\mathcal{D}=|-nK_{V}|$
for $n\gg 0$, $\mathcal{M}=\rho(\mathcal{D})$, and
$M_{X}=\gamma\mathcal{M}$, where $\gamma\in\mathbb{Q}$ such that
$K_{X}+\gamma M_{X}\sim_\mathbb{Q} 0$. Then either $\rho$ is not
biregular or $\mathbb{CS}(X, M_{X})\ne\emptyset$.
\end{theorem}

\begin{proof}
Suppose that $\mathbb{CS}(X, M_{X})=\emptyset$. Then the
singularities of the log pair $(X, M_{X})$ are terminal. In
particular, we have $\kappa(X, M_{X})=0$, which implies
$\gamma={\frac {1}{n}}$. Thus for some rational $\epsilon>\gamma$
the log pair $(X, \epsilon\mathcal{M})$ is canonical model. On the
other hand, the movable log pair $(V, \epsilon\mathcal{D})$ is a
canonical model as well. Hence $\rho$ is biregular by
Lemma~\ref{lemma:canonical-model-is-unique}.
\end{proof}

\section{Birational superrigidity.}
\label{sec:5}

In this section we prove Theorem~\ref{theorem:main}. Let
$\pi:X\to\mathbb{P}^{2n}$ be a cyclic triple cover branched over a
hypersurface $S\subset\mathbb{P}^{2n}$ of degree $3n$ such that
the only singularities of the hypersurface $S$ are ordinary double
points, and $n\ge 2$. Then $X$ is a Fano variety, the
singularities of the variety $X$ are terminal Gorenstein
singularities, and
$K_{X}\sim\pi^{*}(\mathcal{O}_{\mathbb{P}^{2n}}(-1))$.

The variety $X$ is a hypersurface
$$
y^{3}=f_{3n}(x_{0},\ldots,x_{2n})\subset\mathbb{P}(1^{2n+1},n)\cong\mathrm{Proj}(\mathbb{C}[x_{0},\ldots,x_{2n},y]),%
$$
where $f_{3n}$ is a homogeneous polynomial of degree $3n$. The
triple cover $\pi:X\to\mathbb{P}^{2n}$ is a restriction of the
natural projection
$\mathbb{P}(1^{2n+1},n)\dashrightarrow\mathbb{P}^{2n}$ that is
induced by the embedding of the graded algebras
$\mathbb{C}[x_{0},\ldots,x_{2n}]\subset\mathbb{C}[x_{0},\ldots,x_{2n},y]$.
Moreover, the equation of the hypersurface
$S\subset\mathbb{P}^{2n}$ is $f_{3n}(x_{0},\ldots,x_{2n})=0$.

The variety $X$ is $\mathbb{Q}$-factorial, but this must be
proved. We prove a stronger statement following the arguments in
\cite{Ch96a}, \cite{Ch96b}. In fact, the $\mathbb{Q}$-factoriality
of the variety $X$ must follow from the Lefschetz theorem (see
\cite{Bott59}, \cite{AndFr59}, \cite{Gro65}), because $X$ has
isolated singularities.

\begin{lemma}
\label{lemma:factoriality} The groups $\mathrm{Cl}(X)$ and
$\mathrm{Pic}(X)$ are generated by the divisor $K_{X}$.
\end{lemma}

\begin{proof}
Let $D$ be a Weil divisor on $X$. We must show that $D\sim rK_{X}$
for some $r\in\mathbb{Z}$.

Let $H$ be a general divisor in $|-kK_{X}|$ for $k\gg 0$. Then $H$
is a smooth weighted complete intersection in
$\mathbb{P}(1^{2n+1}, n)$ and $\mathrm{dim}(H)\ge 3$. In
particular, $\mathrm{Pic}(H)$ is generated by the divisor
$K_{X}\vert_{H}$ by Theorem~3.13 of chapter~XI in \cite{Gro65}
(see \cite{CalLy94}, Lemma~3.2.2 in \cite{Do82}, or Lemma~3.5 in
\cite{CPR}). Hence there is $r\in\mathbb{Z}$ such that
$D\vert_{H}\sim rK_{X}\vert_{H}$.

Let $\mathrm{\Delta}=D-rK_{X}$. Then the sequence of sheaves
$$
0\to\mathcal{O}_{X}(\mathrm{\Delta})\otimes\mathcal{O}_{X}(-H)\to\mathcal{O}_{X}(\mathrm{\Delta})\to\mathcal{O}_{H}\to 0%
$$
is exact, because $\mathcal{O}_{X}(\mathrm{\Delta})$ is locally
free in the neighborhood of $H$. Thus the sequence
$$
0\to H^{0}(\mathcal{O}_{X}(\mathrm{\Delta}))\to H^{0}(\mathcal{O}_{H})\to H^{1}(\mathcal{O}_{X}(\mathrm{\Delta})\otimes \mathcal{O}_{X}(-H))%
$$
is exact. The sheaf $\mathcal{O}_{X}(\mathrm{\Delta})$ is
reflexive (see \cite{Har80}). So there is an exact sequence of
sheaves
$$
0\to\mathcal{O}_{X}(\mathrm{\Delta})\to\mathcal{E}\to\mathcal{F}\to 0,%
$$
where $\mathcal{E}$ is a locally free sheaf, and $\mathcal{F}$ has
no torsion. Hence the sequence
$$
H^{0}(\mathcal{F}\otimes\mathcal{O}_{X}(-H))\to H^{1}(\mathcal{O}_{X}(\mathrm{\Delta}-H))\to H^{1}(\mathcal{E}\otimes\mathcal{O}_{X}(-H))%
$$
is exact. However, the group
$H^{0}(\mathcal{F}\otimes\mathcal{O}_{X}(-H))$ is trivial because
$\mathcal{F}$ has no torsion, and the group
$H^{1}(\mathcal{E}\otimes\mathcal{O}_{X}(-H))$ is trivial by the
lemma of Enriques--Severi--Zariski (see \cite{Za52}), because the
variety $X$ is normal. Therefore
$$
H^{1}(\mathcal{O}_{X}(\mathrm{\Delta})\otimes\mathcal{O}_{X}(-H))=0%
$$
and $H^{0}(\mathcal{O}_{X}(\mathrm{\Delta}))=\mathbb{C}$.
Similarly, we have
$H^{0}(\mathcal{O}_{X}(-\mathrm{\Delta}))=\mathbb{C}$. Thus $D\sim
rK_{X}$.
\end{proof}

Suppose that the Fano variety $X$ is not birationally superrigid.
Let us show that this assumption leads to a contradiction. It is
follows from Corollary~\ref{corollary:Nother-Fano} that there is a
movable log pair $(X, M_{X})$ such that the set
$\mathbb{CS}(X,M_{X})$ is not empty, the boundary $M_{X}$ is
effective, and the divisor $-(K_{X}+M_{X})$ is ample. The latter
simply means that $M_{X}\sim_{\mathbb{Q}} -rK_{X}$ for a positive
rational number $r<1$. Let $Z\subset X$ be an element of the set
$\mathbb{CS}(X,M_{X})$.

\begin{lemma}
\label{lemma:smooth-points} The subvariety $Z\subset X$ is not a
smooth point of the variety $X$.
\end{lemma}

\begin{proof}
Suppose that $Z$ is a smooth point of $X$. Let
$H_{1},\ldots,H_{2n-2}$ be sufficiently general divisors in the
linear system $|\pi^{*}(\mathcal{O}_{\mathbb{P}^{2n}}(1))|$ such
that each $H_{i}$ passes through $Z$. Then
$$
3>M^{2}_{X}\cdot H_{1}\cdot\cdots \cdot H_{2n-2}\ge \mathrm{mult}_{Z}(M^{2}_{X})\mathrm{mult}_{Z}(H_{1})\cdots\mathrm{mult}_{Z}(H_{2n-2})>4,%
$$
because $\mathrm{mult}_{Z}(M_{X}^{2})>4$ by
Theorem~\ref{theorem:Iskovskikh}, which is a contradiction.
\end{proof}

\begin{lemma}
\label{lemma:singular-points} The subvariety $Z\subset X$ is not a
singular point of the variety $X$.
\end{lemma}

\begin{proof}
Suppose that $Z\subset X$ is a singular point of the variety $X$.
Then $\pi(Z)$ is a singular point of the hypersurface
$S\subset\mathbb{P}^{2n}$. Let $\alpha:V\to X$ be a usual blow up
$Z$ and $G\subset V$ be an $\alpha$-exceptional divisor. Then $V$
is smooth and $G$ is a quadric of dimension $2n-1$ having a single
singular point $O\in G$. Namely, the variety $G\subset V$ is a
quadric cone with the vertex $O\in V$.

Let $M_{V}=\alpha^{-1}(M_{X})$, and $\mathrm{mult}_{Z}(M_{X})$ be
a rational number such that the equivalence
$$
M_{V}\sim_{\mathbb{Q}}
\alpha^{*}(M_{X})-\mathrm{mult}_{Z}(M_{X})G
$$
holds. Then $\mathrm{mult}_{Z}(M_{X})>1$ by
Proposition~\ref{proposition:non-simple-double-point}.

Put $H=\alpha^{*}(-K_{X})$ and consider the linear system $|H-G|$.
By construction, the rational map $\phi_{|H-G|}$ that is given by
the linear system $|H-G|$ is the rational map
$$
\gamma\circ\pi\circ\alpha:V\dasharrow\mathbb{P}^{2n-1},
$$
where $\gamma:\mathbb{P}^{2n}\dasharrow\mathbb{P}^{2n-1}$ is a
projection from the point $\pi(Z)$. The base locus of the linear
system $|H-G|$ is not empty. Namely, its base locus consists of
the vertex $O$ of the quadric cone $G$. Moreover, blowing up the
point $O$ we resolve the indeterminacy of the rational map
$\phi_{|H-G|}$, and the proper transform of the quadric cone $G$
is contracted by $\phi_{|H-G|}$ into the smooth quadric of
dimension $2n-2$.

It should be pointed out that instead of blowing up the points $Z$
and $O$ we can resolve the indeterminacy of the rational map
$\gamma\circ\pi:X\dasharrow\mathbb{P}^{2n-1}$ by a single weighted
blow up
$$
\beta:U\to X\subset\mathbb{P}(1^{2n+1},n)
$$
of the point $Z$ with weights $(2,3^{2n})$ in the corresponding
local coordinates. The weighted blow up $\beta:U\to X$ can be
described as a composition of three rational maps: the blow up
$\alpha$, the blow up of the point $O$, and the consequent
contraction of the proper transform of the quadric cone $G$. The
exceptional divisor of $\beta$ is isomorphic to
$\mathbb{P}^{2n-1}$, and it is a section of the fibration
$\gamma\circ\pi\circ\beta:U\to\mathbb{P}^{2n-1}$. However, the
variety $U$ is singular, namely, the variety $U$ has log terminal
quotient singularities (see \cite{Re87}) of type
${\frac{1}{3}}(1,1)$ along the image of the quadric cone $G$ on
the variety $U$.

Let $C$ be a general curve that is contained in the fibers of
$\phi_{|H-G|}$. Then $C$ is irreducible and reduced, the curve
$\pi\circ\alpha(C)$ is a line passing through the point $\pi(Z)$.
Moreover, we have $C\cdot G=2$, $C\cdot(H-G)=1$, and $O\in C$.
Intersecting the boundary $M_{V}$ with the curve $C$, we obtain
the inequalities
$$
1>3-2\mathrm{mult}_{Z}(M_{X})>M_{V}\cdot C\ge\mathrm{mult}_{O}(M_{V}),%
$$
which imply $\mathrm{mult}_{Z}(M_{X})\le {\frac{3}{2}}$ and
$\mathrm{mult}_{O}(M_{V})<1$. The equivalence
$$
K_{V}+M_{V}\sim_{\mathbb{Q}}\alpha^{*}(K_{X}+M_{X})+(2n-2-\mathrm{mult}_{Z}(M_{X}))G%
$$
and $\mathrm{mult}_{Z}(M_{X})\le {\frac{3}{2}}$ imply the
existence of a proper subvariety $Y\subset G$ such that
$$
Y\in\mathbb{CS}(V, M_{V}-(2n-2-\mathrm{mult}_{Z}(M_{X}))G).
$$

The dimension of $Y$ does not exceed $2n-2$,
$\mathrm{mult}_{Y}(M_{V})>1$, and $Y\in\mathbb{CS}(V, M_{V})$.

Let $\mathrm{dim}(Y)=2n-2$. In the case when $O\in Y$ we have
$$1>\mathrm{mult}_{O}(M_{V})\ge\mathrm{mult}_{Y}(M_{V})>1,$$
which is impossible. Thus $O\not\in Y$. Let $L$ be a general
ruling of the cone $G$. Then
$$
{\frac{3}{2}}\ge\mathrm{mult}_{Z}(M_{X})=M_{V}\cdot L\ge\mathrm{mult}_{Y}(M_{V})L\cdot Y,%
$$
where $L\cdot Y$ is an intersection on $G$. Hence $L\cdot Y=1$ and
$Y$ is a hyperplane section of the quadric cone $G$ under the
natural embedding $G\subset\mathbb{P}^{2n}$. Note, that we have
$$
Y\in\mathbb{LCS}(V, M_{V}-(2n-3-\mathrm{mult}_{Z}(M_{X}))G),
$$
and we can apply Theorem~\ref{theorem:Corti} to the log pair $(V,
M_{V}-(2n-3-\mathrm{mult}_{Z}(M_{X}))G)$ and to the subvariety
$Y\subset V$ of codimension $2$. This gives the inequality
$$
\mathrm{mult}_{Y}(M_{V}^{2})\ge
4(2n-2-\mathrm{mult}_{Z}(M_{X}))\ge 2,
$$
because $\mathrm{mult}_{Z}(M_{X})\le {\frac{3}{2}}$ and $n\ge 2$.
Let $H_{1},\ldots,H_{2n-2}$ be sufficiently general divisors in
the linear system $|H-G|$. Then the inequalities
$$
1>3-2\mathrm{mult}^{2}_{Z}(M_{X})>H_{1}\cdot H_{2}\cdots H_{2n-2}\cdot M^{2}_{V}\ge\mathrm{mult}_{Y}(M^{2}_{V})(H-G)^{2n-2}\cdot Y\ge 2%
$$
hold, which is a contradiction.

Therefore $\mathrm{dim}(Y)<2n-2$. The inequality
$\mathrm{mult}_{O}(M_{V})<1$ implies that $O$ is not contained in
$Y$. Let $P$ be a general point $P\in Y$. Then
$\mathrm{mult}_{P}(M^{2}_{V})>4$ by
Theorem~\ref{theorem:Iskovskikh}.

Let $\mathcal{D}\subset |H-G|$ be a linear subsystem consisting of
divisors that passes through the point $P$. The base locus of the
linear system $\mathcal{D}$ consists of $2$ curves. The first one
is a ruling $L_{P}$ of a quadric cone $G$ such that $P\in L_{P}$.
The second one is a possibly reducible curve $C_{P}$ such that
$\pi\circ\alpha(C_{P})\subset\mathbb{P}^{2n}$ is a line that
passes through the point $\pi(Z)$.

The line $\pi\circ\alpha(C_{P})$ gives a point in the
projectivization of the tangent cone of the hypersurface $S$ at
the point $\pi(Z)$ that corresponds to the image of the point
$\zeta(P)$, where $\zeta$ is a projection of the cone $G$ to its
base. Note that the base of the cone $G$ is canonically isomorphic
to a projectivization of the tangent cone to $S$ at $\pi(Z)$.

Let $D_{1},\ldots,D_{2n-2}$ be general divisors in $\mathcal{D}$,
and $T$ be a one-cycle $H_{1}\cdots H_{2n-3}\cdot M_{V}^{2}$ on
the variety $V$. Then $T$ is an effective and
$\mathrm{mult}_{P}(T)>4$. Unfortunately, we are unable to
intersect properly the cycle $T$ with the remaining divisor
$H_{2n-2}$, because $H_{2n-2}$ may contain components of the
effective one-dimensional cycle $T$. Namely, $H_{2n-2}$ may
contain either the curve $L_{P}$ or components of the the curve
$C_{P}$ in the case if some of them are contained in
$\mathrm{Supp}(T)$.

Suppose that the curve $C_{P}$ is irreducible. Then
$$
T=\mu L_{P} + \lambda C_{P}+\Gamma,
$$
where $\mu$ and $\lambda$ are nonnegative rational numbers, and
$\Gamma$ is an effective one-dimensional cycle whose support does
not contain the curves  $L_{P}$ и $C_{P}$. Then
$$
\mathrm{mult}_{P}(\Gamma)>
4-\mathrm{mult}_{P}(L_{P})\mu-\mathrm{mult}_{P}(C_{P})\lambda=
4-\mu-\mathrm{mult}_{P}(C_{P})\lambda\ge 4-\mu-3\lambda,%
$$
because $\mathrm{mult}_{P}(C_{P})\le 3$, which is implied by the
fact that $C_{P}$ is a triple cover of a line that is blown up in
a possible singular point. Intersecting the effective cycle
$\Gamma$ with the divisor $H_{2n-2}$, we obtain the inequalities
$$
3-2\mathrm{mult}^{2}_{Z}(M_{X})-\mu>\Gamma\cdot H_{2n-2}\ge\mathrm{mult}_{P}(\Gamma)>4-\mu-3\lambda,%
$$
because $C_{P}\cdot H_{2n-2}=0$. Therefore $\lambda>1$.
Intersecting the cycle $T$ with a sufficiently general divisor $H$
of the free linear system $|\alpha^{*}(-K_{X})|$, we immediately
obtain a contradiction, because $H\cdot C_{P}=3$ and $H\cdot T<3$.

Hence the curve $C_{P}$ is reducible. However, the triple cover
$\pi$ is cyclic, which implies
$$
C_{P}=C_{1}+C_{2}+C_{3},
$$
where $C_{i}$ is a nonsingular rational curve such that
$\pi\circ\alpha(C_{P})$ is a line, the restriction morphism
$\pi\circ\alpha\vert_{C_{i}}$ is an isomorphism, $-K_{X}\cdot
\alpha(C_{i})=1$, and $C_{i}\ne C_{j}$ if $i\ne j$. Put
$$
T=\mu L_{P} + \sum_{i=1}^{3}\lambda_{i} C_{i}+\Gamma,
$$
where $\mu$ and $\lambda_{i}$ are nonnegative rational numbers,
and $\Gamma$ is an effective one-dimensional cycle whose support
does not contain the curves $L_{P}$ and $C_{i}$. As in the case of
the irreducible curve $C_{P}$ we can intersect properly the cycle
$\Gamma$ with the divisor $H_{2n-2}$, which immediately implies
the inequality $\sum_{i=1}^{3}\lambda_{i}>1$. Intersecting the
cycle $T$ with a general divisor $H$ in $|\alpha^{*}(-K_{X})|$, we
obtain a contradiction, because $H\cdot C_{i}=1$ and $H\cdot T<3$.
\end{proof}

\begin{lemma}
\label{lemma:codimension-big} The inequality
$\mathrm{codim}(Z\subset X)>2$ is impossible.
\end{lemma}

\begin{proof}
Suppose that $\mathrm{codim}(Z\subset X)>2$. Then
$\mathrm{dim}(Z)\ne 0$ by Lemmas~\ref{lemma:smooth-points} and
\ref{lemma:singular-points}, and
$$
\mathrm{mult}_{Z}(M_{X}^{2})\ge 4
$$
by Theorem~\ref{theorem:Iskovskikh}. Let $O$ be a general point on
$Z$, and let $H_{1},\ldots,H_{n-2}$ be sufficiently general
divisors in $|-K_{X}|$ such that each $H_{i}$ passes through the
point $O$. Then
$$
3>M^{2}_{X}\cdot H_{1}\cdots H_{2n-2}\ge\mathrm{mult}_{Z}(M^{2}_{X})\ge 4,%
$$
which is impossible.
\end{proof}

Therefore we proved that $\mathrm{codim}(Z\subset X)=2$.

\begin{lemma}
\label{lemma:codimension-two-I} The inequality $K_{X}^{2n-2}\cdot
Z\le 2$ holds.
\end{lemma}

\begin{proof}
The inequality $K_{X}^{2n-2}\cdot Z\le 2$ easily follows from the
equality $K_{X}^{2n}=3$, the ampleness of the divisor
$-(K_{X}+M_{X})$, and the inequality $\mathrm{mult}_{Z}(M_{X})\ge
1$.
\end{proof}

\begin{lemma}
\label{lemma:codimension-two-II} The equality $n=2$ holds, namely,
$\mathrm{dim}(X)=4$.
\end{lemma}

\begin{proof}
Suppose that $n>2$. Let $V$ be a general divisor in $|-K_{X}|$.
Then $V$ is a smooth hypersurface of degree $3n$ and of dimension
$2n-1\ge 5$ in $\mathbb{P}(1^{2n},n)$. Hence the cohomology group
$H_{4n-6}(V,\mathbb{C})$ is one-dimensional (see \cite{Ste87},
Theorem~7.2 in \cite{IF00}, and \S 4 in \cite{Do82}).

Let us show that the subvariety
$$
Y=Z\cap V\subset V
$$
of dimension $2n-3$ can not generate the group
$H_{4n-6}(V,\mathbb{C})$. Let $Y\equiv\lambda D^{2}$ in
$H_{4n-6}(V,\mathbb{C})$ for some $\lambda\in\mathbb{C}$, where
$D=-K_{X}\vert_{V}$. The image $\pi(Z)\subset\mathbb{P}^{2n}$ is
either a linear subspace of dimension $2n-2$ or a quadric of
dimension $2n-2$ by Lemma~\ref{lemma:codimension-two-I}.  In
particular, applying the Lefschetz theorem to a smooth hyperplane
section of $S$, we see that $\pi(Z)\not\subset S$.

The subvariety $\pi^{-1}(\pi(Z))$ splits into three irreducible
subvarieties, which are conjugated by the action of the group
$\mathbb{Z}_{3}$ on the variety $X$ that interchanges the fibers
of the triple cover $\pi$. Therefore $\lambda={\frac
{\alpha}{3}}$, where $\alpha=K_{X}^{2n-2}\cdot Z=1$ or $2$ by
Lemma~\ref{lemma:codimension-two-I}. The equality
$$
\alpha=Y\cdot D^{2n-3}=\lambda^{2-n}D\cdot Y^{n-2}
$$
implies
$$
D\cdot Y^{n-2}={\frac {\alpha^{n-1}}{3}}\not\in\mathbb{Z},
$$
but $D\cdot Y^{n-2}\in\mathbb{Z}$, because $V$ is smooth, which is
a contradiction.
\end{proof}

It should be pointed out that we can obtain the claim of
Lemma~\ref{lemma:codimension-two-II} by mean of the applying
Proposition~5 in \cite{Pu02a} or Proposition~4.4 in \cite{dFEM03}
to $S\subset\mathbb{P}^{2n}$ and $S\cap \pi(Z)$.

\begin{lemma}
\label{lemma:codimension-two-III} The surface $\pi(Z)$ is not
contained in the hypersurface $S\subset\mathbb{P}^{4}$.
\end{lemma}

\begin{proof}
In the smooth case the claim is trivial due to the Lefschetz
theorem.

Let $V\subset X$ be a general divisor in the linear system
$|-K_{X}|$. Then the induced morphism
$\tau=\pi\vert_{V}:V\to\mathbb{P}^{3}$ is a cyclic triple cover
branched over a smooth hypersurface
$$
F=S\cap
\pi(V)\subset\mathbb{P}^{3}
$$
of degree $6$. Put
$M_{V}=M_{X}\vert_{V}$ and $C=Z\cap V$. Then the boundary $M_{V}$
is movable and effective, the curve $C$ is smooth and rational,
the curve $\tau(C)$ is either a line or a conic, and the
restriction morphism $\tau\vert_{C}$ is biregular. Moreover, the
inequality $\mathrm{mult}_{C}(M_{V})\ge 1$ and the equivalence
$M_{V}\sim_{\mathbb{Q}} rH$ hold, where $H\sim
\tau^{*}(\mathcal{O}_{\mathbb{P}^{3}}(1))$ and $r\in\mathbb{Q}\cap
(0,1)$.

Suppose that $\tau(C)\subset F$. Let us show that the latter
assumption leads to a contradiction. Let $O$ be a point on $C$.
Put $P=\tau(O)\in \tau(C)$. Let $T\subset\mathbb{P}^{3}$ be a
hyperplane that tangents the hypersurface $F$ at the point $P$.
Then the curve $Y=T\cap F$ is singular at the point $P$. In the
case when the multiplicity of the curve $Y$ in the point $P$ is
$2$, let $L$ be a line in $T$ that passes through the point $P$ in
the direction corresponding to any point in the projectivization
of the tangent cone to the curve $T$ at the point $P$. In the case
when the multiplicity of the curve $Y$ in the point $P$ is greater
than $2$, let $L$ be any line in $T$ that passes through the point
$P$. By construction, the line $L$ tangents $F$ such that the
multiplicity of the tangency is at least $3$.

Let ${\tilde L}=\tau^{-1}(L)$. Then $\mathrm{mult}_{O}({\tilde
L})=3$. Intersecting the curve ${\tilde L}$ with a movable
boundary $M_{V}$, we immediately obtain the following: at least
one of the irreducible components of the curve ${\tilde L}$ is
contained in the base locus of one of the component of the movable
boundary $M_{V}$. However, the latter is impossible in the case
when the line $L$ spans at least a divisor in $\mathbb{P}^{3}$
when we vary the point $O$ on $C$. Let us show that the line $L$
always spans at least a divisor in $\mathbb{P}^{3}$ when we vary
the point $O$ on the curve $C$, which concludes the proof.

It should be pointed out that the hyperplane $T$ tangents the
hypersurface $F$ in finite number of points. This is easily
implied either by the Zak theorem on the finiteness of the Gauss
map (see \cite{FuLa81}, \cite{Ish82}, \cite{Zak93}) or by
Theorem~2 in \cite{Pu95} (see Lemma~3.18 in \cite{Ch03b}).

Suppose that $\tau(C)$ is a line. Then $\tau(C)\subset Y\subset
T$, and $T$ spans a pencil of hyperplanes in $\mathbb{P}^{3}$ that
pass through $\tau(C)$ when we vary $O$ on $C$. Put $Y=\tau(C)\cup
R$. In the case when the point $O$ is sufficiently general, the
curve $R$ is smooth and intersects $\tau(C)$ transversally by the
Bertini theorem. In particular, we always can choose the line $L$
different from the line $\tau(C)$. Therefore different choices of
the sufficiently general point $O$ on $C$ give different line $L$.
Hence the line $L$ spans a divisor when we vary the point $O$ on
$C$.

So $\tau(C)$ is a conic. Then $\tau(C)\not\subset Y$ when $O$ is a
general point on $O$. On the other hand, the hyperplane $T$
tangents $\tau(C)$ in $P$. Therefore, the hyperplane $T$
intersects the conic $\tau(C)$ just by the point $P$ if the point
$O$ is general on $O$. However, the line $L$ passes through the
point $P$, and $L$ is contained in the hyperplane $T$. Thus the
different choices of the sufficiently general point $O$ give
different line $L$.  Hence the line $L$ spans a divisor when we
vary the point $O$ on the curve $C$, which concludes the proof.
\end{proof}

\begin{lemma}
\label{lemma:codimension-two-IV} The surface $\pi(Z)$ is not a
plane in $\mathbb{P}^{4}$.
\end{lemma}

\begin{proof}
Suppose that $\pi(Z)$ is a two-dimensional linear subspace of
$\mathbb{P}^{4}$. Let us consider the reduction to a smooth 3-fold
to get a contradiction as in the proof of
Lemma~\ref{lemma:codimension-two-III}.

Let $V\subset X$ be a general divisor in $|-K_{X}|$, and
$$
\tau=\pi\vert_{V}:V\to\mathbb{P}^{3}
$$
be an induced cyclic triple cover branched over a smooth
hypersurface $F=S\cap \pi(V)\subset\mathbb{P}^{3}$ of degree $6$.
Put $M_{V}=M_{X}\vert_{V}$ and $C=Z\cap V$. Then $M_{V}$ is an
effective movable boundary, the curve $\tau(C)$ is a line, the
morphism $\tau\vert_{C}$ is biregular, the curve $\tau(C)$ is not
contained in the hypersurface $F$, $\mathrm{mult}_{C}(M_{V})\ge
1$, and $M_{V}\sim_{\mathbb{Q}} rH$, where $H\sim
\tau^{*}(\mathcal{O}_{\mathbb{P}^{3}}(1))$ and $r$ is a positive
rational number such that $r<1$. The variety $V$ is a smooth
Calabi-Yau variety, namely, the rational equivalence $K_{V}\sim 0$
holds.

Let $\mathcal{D}\subset
|\tau^{*}(\mathcal{O}_{\mathbb{P}^{3}}(1))|$ be a pencil
consisting of surfaces passing through $C$. The base locus of the
pencil $\mathcal{D}$ consists of the curve $C$ and $2$ different
curves $\tilde{C}$ and $\hat{C}$ such that
$$
\tau(C)=\tau(\tilde{C})=\tau(\hat{C}).
$$

The curves $C$, $\tilde{C}$ and $\hat{C}$ are conjugated by the
action of the group $\mathbb{Z}_{3}$ on $V$ that interchanges the
fibers of the cyclic triple cover $\tau$.

Let $f:U\to V$ be a blow up of $C$ and $E=f^{-1}(C)$. Put
$\mathcal{P}=f^{-1}(\mathcal{D})$. Then
$$
\mathcal{P}\sim D-E,
$$
where $D=(\tau\circ f)^{*}(\mathcal{O}_{\mathbb{P}^{3}}(1))$. On
the other hand, the base locus of the pencil $\mathcal{P}$
consists of proper transforms of the curves $\tilde{C}$ and
$\hat{C}$ on the variety $U$. In particular, the proper transforms
of the curves $\tilde{C}$ and $\hat{C}$ on $U$ are the only curves
on the variety $U$ that has negative intersection with the divisor
$D-E$. Therefore the divisor $2D-E$ is numerically effective on
the variety $U$. In particular, the inequality
$$
(2D-E)\cdot
M_{V}^{2}\ge 0
$$
holds, where $M_{U}$ is a proper transform of the
movable boundary $M_{V}$ on the variety $U$.

Now we calculate the intersection $(2D-E)\cdot M_{U}^{2}\ge 0$
implicitly. Firstly, the equalities
$$
D^{3}=3, D^{2}\cdot E=0, D\cdot E^{2}=-1
$$
hold. Secondly, the equalities
$$E^{3}=-\mathrm{deg}(\mathcal{N}_{C/V})=K_{V}\cdot C+2-2g(C)=2$$
holds (see \cite{Is80a}). Thirdly, the equivalence
$M_{U}\sim_{\mathbb{Q}} rD-\mathrm{mult}_{C}(M_{V})E$ holds. Hence
$$
(2D-E)\cdot M_{U}^{2}=6r^{2}-2\mathrm{mult}^{2}_{C}(M_{V})-2r\mathrm{mult}_{C}(M_{V})-2\mathrm{mult}^{2}_{C}(M_{V}),%
$$
which implies $(2D-E)\cdot M_{U}^{2}<0$, because $r<1$ and
$\mathrm{mult}_{C}(M_{V})\ge 1$.
\end{proof}

\begin{lemma}
\label{lemma:codimension-two-V} The surface $\pi(Z)$ is not a
quadric in $\mathbb{P}^{4}$.
\end{lemma}

\begin{proof}
Suppose that $\pi(Z)$ is an irreducible two-dimensional quadric in
$\mathbb{P}^{4}$. Then we can obtain the contradiction in the same
way as in the proof of Lemma~\ref{lemma:codimension-two-IV}. Let
us show the small modifications that must be done to the proof of
Lemma~\ref{lemma:codimension-two-IV}. We use the notation of the
proof  of Lemma~\ref{lemma:codimension-two-IV}.

Firstly, the curve $\tau(C)$ is a conic. Secondly, the base locus
of the linear system $|2D-E|$ is contained in the union
$\tilde{C}\cup\hat{C}$, because $|2D-E|$ contains proper
transforms of quadric cones in $\mathbb{P}^{3}$ over the conic
$\tau(C)$. However, the intersections of the proper transforms of
the curves $\tilde{C}$ and $\hat{C}$ on $U$ with $2D-E$ are
non-negative. In particular, the divisor $2D-E$ is numerically
effective as in the proof of Lemma~\ref{lemma:codimension-two-IV}.
Thus the inequality $(2D-E)\cdot M_{U}^{2}\ge 0$ holds. Thirdly,
the equality $D\cdot E^{2}=-2$ holds, but $E^{3}=2$. Fourthly, the
equality
$$
(2D-E)\cdot M_{U}^{2}=6r-4\mathrm{mult}^{2}_{C}(M_{V})-4r\mathrm{mult}_{C}(M_{V})-2\mathrm{mult}^{2}_{C}(M_{V})%
$$
holds, which implies $(2D-E)\cdot M_{U}^{2}<0$, because $r<1$ and
$\mathrm{mult}_{C}(M_{V})\ge 1$.
\end{proof}

Therefore, Theorem~\ref{theorem:main} is proved.

\section{The absence of elliptic structures.}
\label{sec:6}

In this section we prove Theorem~\ref{theorem:second}. Let
$\pi:X\to\mathbb{P}^{2n}$ be a cyclic triple cover branched over a
hypersurface $S\subset\mathbb{P}^{2n}$ of degree $3n$ such that
the only singularities of the hypersurface $S$ are ordinary double
points, and $n\ge 2$. Then $X$ is a Fano variety, the
singularities of the variety $X$ are terminal and
$\mathbb{Q}$-factorial (see Lemma~\ref{lemma:factoriality}), and
the equivalence
$$
K_{X}\sim\pi^{*}(\mathcal{O}_{\mathbb{P}^{2n}}(-1))
$$
holds. Suppose that there are birational map
$\rho:\hat{X}\dasharrow X$ and morphism $\nu:\hat{X}\to W$ such
that the generic fiber of $\nu$ is a smooth elliptic curve. Let us
show that the latter assumption leads to a contradiction.

Let $D$ be a very ample divisor $D$ on $W$. Put
$\mathcal{D}=|\nu^{*}(D)|$, $\mathcal{M}=\rho(\mathcal{D})$, and
$M_{X}=\gamma\mathcal{M}$, where $\gamma\in\mathbb{Q}$ such that
$\gamma M_{X}\sim_\mathbb{Q}-K_{X}$. Then $\mathbb{CS}(X,
M_{X})\ne\emptyset$ by
Theorem~\ref{theorem:Nother-Fano-inequality-elliptic-case}.

\begin{remark}
\label{remark:canonical-singularities} It follows from the proof
of Theorem~\ref{theorem:main} that the singularities of the
movable log pair $(X, M_{X})$ are canonical (see
Теорему~\ref{theorem:Nother-Fano-inequality}).
\end{remark}

The claim of Theorem~\ref{theorem:second} is a limit of the claim
of Theorem~\ref{theorem:main}. Therefore we can repeat almost all
steps of the proof of Theorem~\ref{theorem:main} under slightly
weaker conditions. However, we must modify the proof the proof of
Theorem~\ref{theorem:main} using the following property of $(X,
M_{X})$.

\begin{remark}
\label{remark:not-composed-from-a-pencil} The linear system
$\mathcal{M}$ is not composed from a pencil\footnote{Namely, the
inequality $\mathrm{dim}(\psi_{\mathcal{M}}(X))>1$ holds.}.
\end{remark}

Let $Z\subset X$ be an element of the set $\mathbb{CS}(X, M_{X})$.

\begin{proposition}
\label{proposition:compilation} The equality
$\mathrm{codim}(Z\subset X)=2$ holds.
\end{proposition}

\begin{proof}
The claim is implied by the proofs of
Lemmas~\ref{lemma:smooth-points}, \ref{lemma:singular-points},
\ref{lemma:codimension-big}.
\end{proof}

\begin{proposition}
\label{proposition:mulitplicity} The equality
$\mathrm{mult}_{Z}(M_{X})=1$ holds.
\end{proposition}

\begin{proof}
The claim follows from Proposition~\ref{proposition:compilation}
and Remarks~\ref{remark:canonical-singularities} and
\ref{remark:canonical-centers-of-codimension-two}.
\end{proof}

\begin{lemma}
\label{lemma:degree} The inequality $K_{X}^{2n-2}\cdot Z\le 2$
holds.
\end{lemma}

\begin{proof}
The inequality
$$
K_{X}^{2n-2}\cdot Z\le 3
$$ follows from
$K_{X}^{2n}=3$, $M_{X}\sim_{\mathbb{Q}} -K_{X}$ and
$\mathrm{mult}_{Z}(M_{X})=1$.

Suppose that $K_{X}^{2n-2}\cdot Z=3$. Let us show that
$K_{X}^{2n-2}\cdot Z=3$ leads to a contradiction.

Taking an intersection of the cycle $M_{X}^{2}$ with $2n-2$
sufficiently general divisors in the linear system $|-K_{X}|$, we
see that
$$
\mathrm{Supp}(M_{X}^{2})=Z,
$$
where the equality $\mathrm{Supp}(M_{X}^{2})=Z$ does not depend on
the choice of two different divisors in the linear system
$\mathcal{M}$ in the definition of the cycle $M_{X}^{2}$.

Let $P\in X\setminus Z$ be a sufficiently general point, and
$\mathcal{D}\subset\mathcal{M}$ be a linear system consisting of
divisors passing through the point $P$. Then the base locus of the
linear system $\mathcal{D}$ has codimension at least $2$ in $X$,
because $\mathcal{M}$ does not composed from a pencil. Thus
$$
D_{1}\cap D_{2}=Z
$$
in a set-theoretic sense, where $D_{1}$ and $D_{2}$ are
sufficiently general divisors in the linear system $\mathcal{D}$.
Indeed, the divisors $D_{1}$ and $D_{2}$ are contained in the
linear system $\mathcal{M}$ and we have
$\mathrm{Supp}(M_{X}^{2})=Z$. On the other hand, by definition
$P\in D_{1}\cap D_{2}$ and $P\not\in Z$, which is a contradiction.
\end{proof}

It should be pointed out that the proof of
Lemma~\ref{lemma:codimension-two-II} requires that the following
properties of the subvariety $Z$ hold: $\mathrm{codim}(Z\subset
X)=2$ and $K_{X}^{2n-2}\cdot Z\le 2$.

\begin{corollary}
\label{corollary:dimension-four} The equality $n=2$ holds, namely,
we have $\mathrm{dim}(X)=4$.
\end{corollary}

We must reprove Lemmas~\ref{lemma:codimension-two-III} and
\ref{lemma:codimension-two-IV} under the new conditions. We prove
them using the canonicity of the movable log pair $(X, M_{X})$ and
the fact that $\mathcal{M}$ is not composed from a pencil.
However, the proof of Lemma~\ref{lemma:codimension-two-V} is valid
under the new conditions once we have the claims of
Lemmas~\ref{lemma:codimension-two-III} and
\ref{lemma:codimension-two-IV}.

\begin{corollary}
\label{corollary:quadric-case-simple} The case $\pi(Z)\not\subset
S$ and $K_{X}^{2}\cdot Z=2$ is impossible.
\end{corollary}

Hence we must get rid of the following three cases:
\begin{itemize}
\item $\pi(Z)\not\subset S$ and $K_{X}^{2}\cdot Z=1$; %
\item $\pi(Z)\subset S$ and $K_{X}^{2}\cdot Z=1$;%
\item $\pi(Z)\subset S$ and $K_{X}^{2}\cdot Z=2$.
\end{itemize}

\begin{lemma}
\label{lemma:elliptic-case-I} The case $\pi(Z)\not\subset S$ and
$K_{X}^{2}\cdot Z=1$ is impossible.
\end{lemma}

\begin{proof}
Suppose that $\pi(Z)\not\subset S$ and $K_{X}^{2}\cdot Z=1$. Let
us show that this assumption leads to a contradiction. The surface
$\pi(Z)$ is a two-dimensional linear subspace of $\mathbb{P}^{4}$,
which is not contained in the hypersurface $S$. The triple cover
$\pi$ is cyclic, which implies the existence of two different
surfaces $\tilde{Z}$ and $\hat{Z}$ such that
$$
\pi(Z)=\pi(\tilde{Z})=\pi(\hat{Z}),
$$
and three surfaces $Z$, $\tilde{Z}$ and $\hat{Z}$ are conjugate
under the action of the group $\mathbb{Z}_{3}$ on the variety $X$,
that interchanges the fibers of $\pi$.

Let $V\subset X$ be a sufficiently general divisor in the linear
system $|-K_{X}|$ and
$$
\tau=\pi\vert_{V}:V\to\mathbb{P}^{3}
$$ be an
induced cyclic triple cover. Then the triple cover $\tau$ is
branched over a smooth hypersurface $F=S\cap
\pi(V)\subset\mathbb{P}^{3}$ of degree $6$.

Let $\mathcal{H}=\mathcal{M}\vert_{V}$ and
$M_{V}=M_{X}\vert_{V}=\gamma\mathcal{H}$. Then the base locus of
the linear system $\mathcal{H}$ has codimension at least $2$ in
$V$, the equivalence
$$
M_{V}\sim_{\mathbb{Q}} \tau^{*}(\mathcal{O}_{\mathbb{P}^{3}}(1))
$$
holds. Moreover, the generality in the choice of $V$ implies that
$\mathcal{H}$ is not composed from a pencil. Let $C=Z\cap V$,
$\tilde{C}=\tilde{Z}\cap V$, and $\hat{C}=\hat{Z}\cap V$. Then
$\mathrm{mult}_{C}(M_{V})=1$.

Let $f:U\to V$ be a blow up of a smooth curve $C$, and $E$ be an
exceptional divisor of the blow up $f$. Put
$\mathcal{D}=f^{-1}(\mathcal{H})$ and
$M_{U}=f^{-1}(M_{V})=\gamma\mathcal{D}$. Then
$$M_{U}\sim_{\mathbb{Q}} D-E,$$
where $D=(\tau\circ f)^{*}(\mathcal{O}_{\mathbb{P}^{3}}(1))$.
However, the base locus of the pencil $|D-E|$ consists of proper
transforms of the curves $\tilde{C}$ and $\hat{C}$ on the variety
$U$. Moreover, the equalities
$$(D-E)\cdot\tilde{C}=(D-E)\cdot\hat{C}=-1$$
holds. Therefore the proper transforms of the curves $\tilde{C}$
and $\hat{C}$ on the variety $U$ are the only curves on $U$ that
have non-positive intersection with $2D-E$. In particular, the
divisor $2D-E$ is numerically effective and the inequality
$(2D-E)\cdot M_{U}^{2}\ge 0$ holds.

The intersection $(2D-E)\cdot M_{U}^{2}$ can be easily calculated
(see the proof of Lemma~\ref{lemma:codimension-two-IV}), namely,
the equalities
$$
(2D-E)\cdot M_{U}^{2}=6-2\mathrm{mult}^{2}_{C}(M_{V})-2\mathrm{mult}_{C}(M_{V})-2\mathrm{mult}^{2}_{C}(M_{V})=0%
$$
hold. Thus $\mathrm{Supp}(M_{U}^{2})$ is contained in the curves
$\tilde{C}$ and $\hat{C}$. This is simply means that for any two
different divisors $H_{1}$ and $H_{2}$ in the linear system
$\mathcal{D}$, the intersection $H_{1}\cap H_{2}$ is contained in
the union $\tilde{C}\cup\hat{C}$ in a set-theoretic sense.

Let $P\in U\setminus (\tilde{C}\cup\hat{C})$ be a sufficiently
general point, and $\mathcal{P}\subset \mathcal{D}$ be a linear
subsystem of divisors passing through the point $P$. Then the
linear system $\mathcal{P}$ does not have base components, because
$\mathcal{D}$ is not composed from a pencil. Let $D_{1}$ and
$D_{2}$ be two sufficiently general divisors in the linear system
$\mathcal{P}$. Then in a set-theoretic sense
$$P\in D_{1}\cap D_{1}\subset \tilde{C}\cup\hat{C},$$
because $D_{i}\in \mathcal{D}$. The obtained  contradiction
concludes the proof.
\end{proof}

\begin{lemma}
\label{lemma:elliptic-case-II} The case $\pi(Z)\subset S$ and
$K_{X}^{2}\cdot Z=1$ is impossible.
\end{lemma}

\begin{proof}
Suppose that $\pi(Z)\subset S$ and $K_{X}^{2}\cdot Z=1$. Then
$\pi(Z)$ is a two-dimensional linear subspace of $\mathbb{P}^{4}$,
which is contained in the hypersurface $S$. The Lefschetz theorem
implies that the hypersurface $S$ is singular.

We use the reduction to a smooth 3-fold as in the proof of
Lemma~\ref{lemma:elliptic-case-I}. Moreover, let us use the
notations of the proof of Lemma~\ref{lemma:elliptic-case-I}, which
can be used in this case with the only difference that the
surfaces $Z$, $\tilde{Z}$, $\hat{Z}$ are coincide under the new
conditions, because the surface $Z$ is invariant under the action
of $\mathbb{Z}_{3}$ on the variety $X$ that interchanges the
fibers of $\pi$.

All steps of the proof of Lemma~\ref{lemma:elliptic-case-I}
remains valid under new conditions except the very last one.
Namely, the numerical effectivity of the divisor  $2D-E$ is not
clear, but it can be proved analyzing the class of the divisor
$E\vert_{E}$ in the Picard group of the $f$-exceptional surface
$E\cong\mathbb{F}_{k}$. However, we prove the numerical
effectivity of the divisor  $2D-E$ using more geometric ideas.

Let us consider the pencil $|D-E|$ on the variety $U$. The base
locus of $|D-E|$ consists of a curve $\bar{C}\subset E$ such that
$\bar{C}$ is a section of the projection $f\vert_{E}:E\to C$. It
should be pointed out that the curve $\bar{C}\subset E$ is an
infinitesimal analog of the curve $\tilde{C}$ in the proof of
Lemma~\ref{lemma:elliptic-case-I}. Moreover, blowing up the curve
$\bar{C}$, we can obtain an infinitesimal analog of the curve
$\hat{C}$ in the proof of Lemma~\ref{lemma:elliptic-case-I}, but
we do not need this.

Let $Y$ be a general surface $Y$ in the pencil $|D-E|$. Then $Y$
is singular. Let us describe the singularities of the surface $Y$.
The surface
$$
\tau\circ f(Y)\subset \mathbb{P}^{3}
$$
is a plane passing through the line $\tau(C)\subset F$, where $F$
is a ramification surface of the cyclic triple cover
$\tau:V\to\mathbb{P}^{3}$. In particular, the curve
$$
\tau\circ
f(Y)\cap F
$$
is reducible, it consists of two irreducible components: the line
$\tau(C)$ and a plane quintic curve $R$. Moreover, the quintic $R$
is smooth by the Bertini theorem, and $R$ intersects the line
$\tau(C)$ transversally in $5$ points. On the other hand, the
morphism
$$
\tau\vert_{f(Y)}:f(Y)\longrightarrow\tau\circ
f(Y)\cong\mathbb{P}^{2}
$$
is a cyclic triple cover branched over a curve $\tau(C)\cup R$.
Therefore the singularities of the surface $f(Y)$ are $5$ singular
points of type $\mathbb{A}_{2}$ contained in  $C$. The birational
morphism
$$
f\vert_{Y}:Y\longrightarrow f(Y)
$$
partially resolves the singularities of $f(Y)$. Namely, the
surface $Y$ has $5$ ordinary double points, and each of them
dominates the corresponding singular point of the surface $f(Y)$.

Let $M_{Y}=M_{U}\vert_{Y}$. Then the boundary $M_{Y}$ may not be
movable, because it may contain a multiple of the curve $\bar{C}$
as a fixed component. Hence we can put
$$
M_{Y}=\alpha\bar{C}+\Gamma,
$$
where $\alpha\in\mathbb{Q}_{>0}$ and $\Gamma$ is a movable
boundary on the surface $Y$. On the other hand, we have
$M_{Y}\sim_{\mathbb{Q}} 2\bar{C}$. Moreover, it follows from the
subadjunction formula (see \cite{Ko91}) that
$$
\bar{C}^{2}=-3+\deg(\mathrm{Diff}_{\bar{C}}(0))=-3+5{\frac{1}{2}}<0,
$$
which implies $\alpha=2$ and $\Gamma=\emptyset$. So for any
general divisors $D\in\mathcal{D}$ and $H\in |D-E|$, the
intersection $D\cap H$ is contained in the curve $\bar{C}$ in a
set-theoretic sense.

It should be pointed out that we used the following properties of
the boundary $M_{X}$ in the above arguments: $\mathcal{M}$ does
not have fixed components, the equivalence $M_{X}\sim_{\mathbb{Q}}
-K_{X}$ holds, and $\mathrm{mult}_{Z}(M_{X})=1$. In particular, we
did not use the fact that $\mathcal{M}$ is not composed from a
pencil. So we can repeat the previous arguments to any linear
subsystem $\mathcal{B}\subset\mathcal{M}$ such that $\mathcal{B}$
does not have fixed components. Indeed, the equivalence
$\gamma\mathcal{B}\sim_{\mathbb{Q}} -K_{X}$ and the inequality
$\mathrm{mult}_{Z}(\gamma\mathcal{B})\ge 1$ are obvious, and the
proof of Theorem~\ref{theorem:main} implies the canonicity of the
log pair $(X, \gamma\mathcal{B})$, which implies the equality
$\mathrm{mult}_{Z}(\gamma\mathcal{B})=1$. Therefore for any
sufficiently general divisor $B$ in any linear system
$\mathcal{B}\subset\mathcal{M}$ with no fixed components and a
sufficiently general divisor $H\in |D-E|$, the intersection $B\cap
f(H)$ is contained in the curve $C$ in a set-theoretic sense.

Let $P$ be a sufficiently general point in $X\setminus C$, and
$\mathcal{B}\subset \mathcal{D}$ be a linear subsystem of divisors
passing through the point $P$. Then $\mathcal{B}$ does not have
fixed components, because the linear system $\mathcal{M}$ is not
composed from a pencil. Let $B$ be a general divisor in the linear
system $\mathcal{B}$, and $H$ be a general divisor in $|D-E|$.
Then the intersection $B\cap f(H)$ contains the point $P\not\in
C$. Thus $B\cap f(H)$ is not contained in the curve $C$ in a
set-theoretic sense, which is a contradiction.
\end{proof}

\begin{lemma}
\label{lemma:elliptic-case-III} The case $\pi(Z)\subset S$ and
$K_{X}^{2}\cdot Z=2$ is impossible.
\end{lemma}

\begin{proof}
Suppose that $\pi(Z)\subset S$ and $K_{X}^{2}\cdot Z=2$. Let us
show that this assumption leads to a contradiction. The surface
$\pi(Z)$ is a two-dimensional quadric in $\mathbb{P}^{4}$, which
is contained in the sextic $S$, which implies that $Z$ is
invariant under the action of the group $\mathbb{Z}_{3}$ on the
variety $X$ that interchanges the fibers of $\pi$, because
$\pi(Z)\subset S$. The Lefschetz theorem implies that the
hypersurface $S$ is singular.

We reduce the problem to a smooth 3-fold. Let $V\subset X$ be a
sufficiently general divisor in the linear system $|-K_{X}|$, and
$\tau=\pi\vert_{V}:V\to\mathbb{P}^{3}$ be the induced cyclic
triple cover branched over the smooth hypersurface
$$
F=S\cap \pi(V)\subset\mathbb{P}^{3}
$$
of degree $6$. Put $M_{V}=M_{X}\vert_{V}$. Then $M_{V}$ is a
movable boundary, the equivalence
$$
M_{V}\sim_{\mathbb{Q}} \tau^{*}(\mathcal{O}_{\mathbb{P}^{3}}(1))
$$
holds, and $\mathrm{mult}_{C}(M_{V})=1$, where $C=Z\cap V$. The
curve $\tau(C)\subset F$ is a smooth conic.

Let $f:U\to V$ be a blow up of $C$, and $E=f^{-1}(C)$. Put
$M_{U}=f^{-1}(M_{V})$. Then
$$
M_{U}\sim_{\mathbb{Q}} D-E,
$$
where $D=(\tau\circ f)^{*}(\mathcal{O}_{\mathbb{P}^{3}}(1))$. In
the case when the divisor $2D-E$ is numerically effective, the
explicit calculation of the intersection $(2D-E)\cdot M_{U}^{2}\ge
0$ leads to a contradiction in the same way as in the proof of
Lemma~\ref{lemma:codimension-two-V}. So we may assume that $2D-E$
is not numerically effective

The base locus of $|2D-E|$ is contained in $E$, because $|2D-E|$
contains proper transforms of quadric cones in $\mathbb{P}^{3}$
over the conic $C$. Therefore $2D-E$ is numerically effective if
and only if it has non-negative intersection with the exceptional
section of the ruled surface $E$.

Let us show that $(2D-E)\cdot s_{\infty}\ge 0$, where $s_{\infty}$
is an exceptional section of the ruled surface
$E\cong\mathbb{F}_{k}$, which concludes the proof. The curve $C$
is smooth and $C\cong\mathbb{P}^{1}$. Hence
$$
\mathcal{N}_{C/V}\cong\mathcal{O}_{\mathbb{P}^{1}}(a)\oplus\mathcal{O}_{\mathbb{P}^{1}}(b)
$$
for integers $a$ and $b$ such that $b\ge a$. Then $k=b-a$ and the
equalities
$$
a+b=\mathrm{deg}(\mathcal{N}_{C/V})=2g(C)-2-K_{V}\cdot C=-2
$$
and $E^{3}=-\mathrm{deg}(\mathcal{N}_{C/V})=2$ holds. On the other
hand, the smooth curve $C$ is contained in the smooth surface
$\bar{F}=\tau^{-1}(F)$. Thus the sequence of sheaves
$$
0\to\mathcal{N}_{C/\bar{F}}\to\mathcal{N}_{C/V}\to\mathcal{N}_{\bar{F}/V}\to
0
$$
is exact, where
$\mathcal{N}_{C/\bar{F}}\cong\mathcal{O}_{\mathbb{P}^{1}}(-6)$,
because $C^{2}=-6$ on the surface $\bar{F}\cong F$ by the
adjunction formula. Hence $a\ge -6$. Let $l\subset E$ be a fiber
of the projection $f\vert_{E}$. Then
$$
-E\vert_{E}\sim s_{\infty}+rl
$$
for $r={\frac{2+k}{2}}$, because the equalities
$$
2=E^{3}=(s_{\infty}+rl)^{2}=-k+2r,
$$
holds. So we have
$$
(2D-E)\cdot s_{\infty}=4-E\cdot
s_{\infty}=4+(s_{\infty}+{\frac{2+k}{2}}l)\cdot
s_{\infty}=4-k+{\frac{2+k}{2}}={\frac{10-k}{2}}=6+a\ge 0,
$$
which concludes the proof.
\end{proof}

Therefore Theorem~\ref{theorem:second} is proved.

\section{The proof of Theorems~\ref{theorem:third} and \ref{theorem:forth}.}
\label{sec:7}

In this section we prove Theorems~\ref{theorem:third} and
\ref{theorem:forth}. Let $\pi:X\to\mathbb{P}^{2n}$ be a cyclic
triple cover branched over a hypersurface
$S\subset\mathbb{P}^{2n}$ of degree $3n$ such that $n\ge 2$, and
the only singularities of $S$ are isolated ordinary double and
triple points. Namely, the projectivization of the tangent cone to
the hypersurface $S$ at any singular point $P$ of $S$ is a smooth
hypersurface in $\mathbb{P}^{2n-1}$ of degree
$\mathrm{mult}_{P}(S)\le 3$.

\begin{remark}
\label{remark:factoriality} The proof of
Lemma~\ref{lemma:factoriality} implies that the groups
$\mathrm{Pic}(X)$ and $\mathrm{Cl}(X)$ are generated by the
divisor $-K_{X}$, because the singularities of $X$ are isolated.
\end{remark}

Hence $X$ is a Fano variety with terminal $\mathbb{Q}$-factorial
singularities. We must prove the following three results:
\begin{itemize}
\item the variety $X$ is birationally superrigid;%
\item the variety $X$ is not birationally equivalent to any Fano
variety\hfill\break with
canonical singularities that is not biregular to $X$;%
\item if the variety $X$ is birational to an elliptic fibration,
then $n=2$, \hfill\break the hypersurface $S$ has a triple point
$O$ such that the elliptic fibration\hfill\break is induced by the
projection $\gamma:\mathbb{P}^{4}\dasharrow\mathbb{P}^{3}$ from
the point $O$.
\end{itemize}

Suppose that at least one of the above three claims is not true.
Then there is a linear system $\mathcal{M}$ on the variety $X$
that satisfies the following properties:
\begin{itemize}
\item the linear system $\mathcal{M}$ does not have fixed components;%
\item the set $\mathbb{CS}(X, {\frac{1}{d}}\mathcal{M})$ is not
empty, where $d\in\mathbb{N}$ such that $\mathcal{M}\sim -dK_{X}$;%
\item the linear system $\mathcal{M}$ is not composed from a pencil;%
\item in the case when $n=2$, for any point $O\in S$ such that
$\mathrm{mult}_{O}(S)=3$,\hfill\break the linear system
$\mathcal{M}$ is not contained in the fibers of the rational map
$\gamma\circ\pi$,\hfill\break where
$\gamma:\mathbb{P}^{4}\dasharrow\mathbb{P}^{3}$ is a projection
from the point $O$.
\end{itemize}

The existence of $\mathcal{M}$ follows from
Theorems~\ref{theorem:Nother-Fano-inequality-elliptic-case} and
\ref{theorem:Nother-Fano-inequality-Fano-case} and the proof of
Theorem~\ref{theorem:Nother-Fano-inequality}.

Let us show that the linear system $\mathcal{M}$ does not exist.
Let $Z\subset X$ be a subvariety such that $Z\in\mathbb{CS}(X,
{\frac{1}{d}}\mathcal{M})$. Then the proof of
Theorems~\ref{theorem:main} and \ref{theorem:second} implies that
$Z$ is a singular point of the variety $X$ such that $O=\pi(Z)$ is
an ordinary triple point of $S$.

\begin{remark}
\label{remark:tripple-point} The point $Z$ is an ordinary triple
point of the variety $X$.
\end{remark}

Let $\alpha:V\to X$ be a blow up of the point $O$, and
$E=\alpha^{-1}(O)$. Then $E$ is a smooth hypersurface of degree
$3$ in $\mathbb{P}^{2n}$, and $E\vert_{E}\sim H$, where $H$ is a
hyperplane section of the hypersurface $E\subset\mathbb{P}^{2n}$.
Moreover, the linear system
$$
|\alpha^{*}(-K_{X})-E|
$$
is free and gives a morphism  $\psi:V\to\mathbb{P}^{2n-1}$ such
that $\psi=\gamma\circ\pi\circ\alpha$, where
$\gamma:\mathbb{P}^{2n}\dasharrow\mathbb{P}^{2n-1}$ is a
projection from the point $O$. Let
$\mathrm{mult}_{Z}(\mathcal{M})$ be an integer number such that
$$
\mathcal{D}\sim
\alpha^{*}(-dK_{X})-\mathrm{mult}_{Z}(\mathcal{M})E,
$$
where $\mathcal{D}$ is a proper transform of the linear system
$\mathcal{M}$ on the variety $V$. Let $C\subset V$ be a
sufficiently general curve in a fiber of $\psi$. Then
$$
\mathcal{D}\cdot C=3(d-\mathrm{mult}_{Z}(\mathcal{M}))\ge 0,
$$
and the equality $\mathcal{D}\cdot C=0$ implies that $\mathcal{D}$
is contained in the fibers of $\psi$. On the other hand, the
inequality $\mathrm{mult}_{Z}(\mathcal{M})>d$ holds when $n>2$ and
the inequality $\mathrm{mult}_{Z}(\mathcal{M})\ge d$ holds when
$n=2$ by Theorem~\ref{theorem:tripple-point}. Hence $n=2$ and the
linear system $\mathcal{M}$ is contained in the fibers of the
rational map $\gamma\circ\pi$, which contradicts to one of the
properties of the linear system $\mathcal{M}$. Thus both
Theorems~\ref{theorem:third} and \ref{theorem:forth} are proved.

\section{Potential density.}
\label{sec:8}

In this section we prove Theorem~\ref{theorem:five}. Let
$\pi:X\to\mathbb{P}^{4}$ be a cyclic triple cover branched over a
hypersurface $S\subset\mathbb{P}^{4}$ of degree $6$ such that the
hypersurface $S$ is defined over a number field $\mathbb{F}$.
Suppose that the hypersurface $S$ has an ordinary triple point
$O$, and the hypersurface $S$ is smooth outside of the point $O$.
Thus the equality $\mathrm{mult}_{O}(S)=3$ holds, and the
projectivization of the tangent cone to the hypersurface $S$ at
the point $O$ is a smooth cubic surface in $\mathbb{P}^{3}$. The
point $O$ is defined over the field $\mathbb{F}$.

The variety $X$ can be considered as a hypersurface
$$
y^{3}=x_{0}^{3}f_{3}(x_{1},\ldots,x_{4})+x_{0}^{2}f_{4}(x_{1},\ldots,x_{4})+x_{0}f_{5}(x_{1},\ldots,x_{4})+f_{6}(x_{1},\ldots,x_{4})
$$
in the weighted projective space
$\mathbb{P}(1^{5},2)\cong\mathrm{Proj}(\mathbb{F}[x_{0},\ldots,x_{4},y])$,
where $f_{i}$ is a homogeneous polynomial of degree $i$. The
cyclic triple cover $\pi:X\to\mathbb{P}^{4}$ is a restriction to
the hypersurface $X$ of the natural projection
$\mathbb{P}(1^{5},2)\dashrightarrow\mathbb{P}^{4}$ that is induced
by the natural embedding of the graded algebras
$\mathbb{F}[x_{0},\ldots,x_{4}]\subset\mathbb{F}[x_{0},\ldots,x_{4},y]$.
Moreover, the hypersurface $S\subset\mathbb{P}^{4}$ is given by
the equation
$$
x_{0}^{3}f_{3}(x_{1},\ldots,x_{4})+x_{0}^{2}f_{4}(x_{1},\ldots,x_{4})+x_{0}f_{5}(x_{1},\ldots,x_{4})+f_{6}(x_{1},\ldots,x_{4})=0,
$$
where the coordinates of the singular point $O$ are
$(1:0:\cdots:0)$.

\begin{remark}
\label{remark:tangent-cone} The equation
$f_{3}(x_{1},\ldots,x_{4})=0$ defines a smooth cubic surface in
$\mathbb{P}^{3}$, which is a projectivization of the tangent cone
to $S$ at the point $O$. In particular,  $f_{3}$ is irreducible.
\end{remark}

Suppose that $X$ satisfies the following generality conditions:
\begin{enumerate}
\item $f_{4}$ is not divisible by $f_{3}$;%
\item $f_{5}^{2}-3f_{4}f_{6}$ and
$f_{4}^{2}f_{5}^{2}-4f_{4}^{3}f_{6}-4f_{3}f_{5}^{3}+18f_{3}f_{4}f_{5}f_{6}-27f_{3}^{2}f_{6}^{2}$
are coprime.
\end{enumerate}

\begin{remark}
\label{remark:generality-condition-is-OK} The required generality
conditions are satisfied in the case when the polynomial $f_{i}$
are chosen sufficiently general. The geometrical meaning of the
generality conditions are the following:
\begin{enumerate}
\item a sufficiently general line $L$ in $\mathbb{P}^{4}$ that
passes through the point $O$ and that is contained in the tangent
cone to the hypersurface $S$ at the point $O$ intersects the
hypersurface $S$ in two points that are different from $O$;%
\item there is at most one-dimensional family of curves $C\subset
X$ such that the singular point $P=\pi^{-1}(O)$ of the variety $X$
is contained in the curve
$C$ and $-K_{X}\cdot C=1$.%
\end{enumerate}
\end{remark}

We use the methods of \cite{BoTsch99}, \cite{HaTsch00},
\cite{BoTsch00} to prove the following result implying
Theorem~\ref{theorem:five}.

\begin{proposition}
\label{proposition:explicit-potential-density} Under the
generality conditions, the rational points on $X$ are potentially
dense\footnote{To be precise we must say that $\mathbb{F}$-points
are potentially dense on the variety $X$.}, namely, the set of all
$\mathbb{K}$-points of the variety $X$ is Zariski dense in $X$ for
a finite extension of fields $\mathbb{F}\subset\mathbb{K}$.
\end{proposition}

There are two ways of looking at the potential density of rational
points. The optimistic point of view is the following: the
potential density of rational points reflects the measure of how
close a given variety to being rational. For example, the
geometrical rationality obviously implies the potential density of
rational points. From this point of view the claim of
Theorem~\ref{theorem:five} is very natural, as well as the fact
that we are unable to prove the potential density of rational
points on many many rationally connected nonrational varieties.
For example, it is unknown whether rational points are potentially
dense on the generic quintic hypersurface in $\mathbb{P}^{5}$ or
not (see \cite{Pu87}, \cite{Ch00b}, \cite{dFEM03}). The
pessimistic point of view considers the potential density of
rational points as a much weaker birational invariant. In
particular, there is the following conjecture (see
\cite{HaTsch00}).

\begin{conjecture}
\label{conjecture:potential-density} Let $V$ be a smooth variety
such that $V$ is defined over a number field, and $-K_{V}$ is
numerically effective. Then rational points on $V$ are potentially
dense.
\end{conjecture}

Therefore from the point of view of
Conjecture~\ref{conjecture:potential-density} the claim of
Proposition~\ref{proposition:explicit-potential-density} is just
an illustration of a general principle. It is known that
Conjecture~\ref{conjecture:potential-density} holds the following
algebraic varieties: abelian varieties (see \cite{Has03}), smooth
Fano 3-folds except a double cover of $\mathbb{P}^{3}$ ramified in
a smooth sextic (see \cite{BoTsch99}, \cite{HaTsch00}), smooth
Enriques surfaces (see \cite{BoTsch98}), smooth elliptic K3
surfaces (see \cite{BoTsch00}), smooth K3 surfaces with an
infinite group of auto\-mor\-phisms (see \cite{BoTsch00}), some
symmetric products (see \cite{HasTsch99}). Therefore rational
points are potentially dense on many varieties that are not
rationally connected. However, it is unknown whether rational
points are potentially dense on a generic double cover of
$\mathbb{P}^{2}$ branched over a smooth quartic curve or not (see
\cite{BoTsch99}).

\begin{example}
\label{example:Faltings} Let $C$ be a smooth connected curve such
that the curve $C$ is defined over a number field, and $g(C)\ge
2$. Then rational points are not potentially dense on
$C\times\mathbb{P}^{k}$ by the Faltings theorem (see \cite{Fa83},
\cite{FaWu84}).
\end{example}

It is natural to expect that the potential density of rational
points reflects such birational properties of an algebraic variety
as rational connectedness. However, even in the case of a smooth
conic bundle $\zeta:V\to\mathbb{P}^{n}$ with sufficiently general
and big discriminant it is not known whether rational points are
potentially dense on $V$ or not in the case $n\ge 2$, but it is
known that the potential density of rational points on $V$ is
implied by the Schinzel conjecture for $\zeta:V\to\mathbb{P}^{n}$
(see \cite{CoSw94}). The variety $V$ is nonrational (see
\cite{Sa80} и \cite{Sa82}) and it is expected that $V$ is not
unirational. Perhaps, the potential density of rational points can
be used to obtain an example of a rationally connected variety
that is not unirational.

An example in \cite{CoSkSw97} implies the following generalization
of Conjecture~\ref{conjecture:potential-density}.

\begin{conjecture}
\label{conjecture:potential-density-big} Let $V$ be a smooth
variety such that $V$ is defined over a number field, and the
divisor $-K_{V}$ is not numerically effective. Then rational
points on $V$ are potentially dense if there is no unramified
finite morphism  $f:U\to V$ such that there is a dominant rational
map $g:U\dasharrow Z$, where $Z$ is a variety of general type of
dimension $\mathrm{dim}(Z)>0$.
\end{conjecture}

It should be pointed out that both
Conjectures~\ref{conjecture:potential-density} and
\ref{conjecture:potential-density-big} are logical negation of the
following weak Lang conjecture, which is proved only for curves
and subvarieties of abelian varieties (see \cite{Fa83},
\cite{FaWu84}, \cite{Fa91}).

\begin{conjecture}
\label{conjecture:Lang} Let $V$ be a smooth variety of general
type such that the variety $V$ is defined over a number field.
Then rational points on $V$ are not potentially dense.
\end{conjecture}

The claim of Theorem~\ref{theorem:five} must remain valid without
any generality conditions. Moreover, in non-general case the proof
of the potential density of rational points must be easier than in
general case. The same can be said about the singularities.
Namely, the proof of the potential density of rational points must
become easier when the singularities become worse. However, there
are exceptional extreme cases.

\begin{example}
\label{example:degenerate-triple-cover} Let
$\chi:Y\to\mathbb{P}^{4}$ be a cyclic triple cover branched over a
hypersurface $G$ of degree $6$ such that $G$ is a union of $6$
different hyperplanes defined over a number field $\mathbb{F}$ and
passing through some two-dimensional linear subspace
$\Pi\subset\mathbb{P}^{4}$. Then $Y$ is birational to the product
$C\times \mathbb{P}^{3}$, where $C$ is a cyclic triple cover of
$\mathbb{P}^{1}$ branched over $6$ points, which are defined over
the field $\mathbb{F}$. Then rational points on the variety $Y$
are not potentially dense, because $g(C)=4$ (see
Example~\ref{example:Faltings}).
\end{example}

Let us prove
Proposition~\ref{proposition:explicit-potential-density}. The
following result is due to \cite{Mer96}.

\begin{theorem}
\label{theorem:Merel} Let $\mathbb{F}$ be a number field. Then
there is $n(\mathbb{F})\in\mathbb{N}$ such that $n(\mathbb{F})$
depends only on the field $\mathbb{F}$, and the order of any
torsion $\mathbb{F}$-point on any elliptic curve $C$ is bounded by
$n(\mathbb{F})$, where $C$ is defined over $\mathbb{F}$.
\end{theorem}

Let $P=\pi^{-1}(O)$. Then $P$ is an ordinary triple point on $X$.
Let $\alpha:U\to V$ be a blow up of $P$, and $E$ be an exceptional
divisor of $\alpha$. Then $-K_{U}\sim \alpha^{*}(-K_{X})-E$, the
linear system $|-K_{U}|$ has no base points. Let
$$
\psi:U\to\mathbb{P}^{3}
$$
be a morphism that is given by $|-K_{U}|$. Then $\psi$ is an
elliptic fibration such that $E$ is a three-section of $\psi$, and
$\psi=\gamma\circ\pi$, where
$\gamma:\mathbb{P}^{4}\dasharrow\mathbb{P}^{3}$ is a projection
from the point $O$.

\begin{remark}
\label{remark:unirationality-of-cubic-threefold} The variety $E$
is a smooth cubic hypersurface in $\mathbb{P}^{4}$. The cubic $E$
is not rational over $\mathbb{C}$ (see \cite{ClGr72}), but $E$ is
unirational over $\mathbb{C}$ (see \cite{Ma72}). In particular,
rational points on the variety $E$ are potentially dense.
\end{remark}

Let $D$ be an intersection of two general divisors in $|-K_{U}|$.
Then $D$ is a smooth elliptic surface. The restriction
$\tau=\psi\vert_{D}:D\to \mathbb{P}^{1}$ is a canonical morphism
of the suraface $D$, namely, the equivalence
$K_{U}\sim\tau^{*}(\mathcal{O}_{\mathbb{P}^{1}}(1))$ holds. The
curve $Z=E\cap D$ is a smooth elliptic curve, and the restriction
$\tau\vert_{Z}:Z\to\mathbb{P}^{1}$ is a cyclic triple cover
branched over three points.

\begin{remark}
\label{remark:no-reducible-fibers} The proper transform on the
variety $V$ of every irreducible component of any reducible fiber
of the fibration $\tau$ is a rational curve whose intersection
with the anticanonical divisor $-K_{X}$ is equal to $1$. The
generality conditions implies that there is no more than
one-dimensional family of such curves on $V$. On the other hand,
the generality in the choice of the surface $D$ in the fibers of
$\psi$ and the equality $\mathrm{codim}(D\subset U)=2$ imply that
all fibers of the fibration $\tau$ are irreducible.
\end{remark}

Let $F_{1}$, $F_{2}$, $F_{3}$ be fibers of $\tau$ that pass
through the ramification points of the triple cover
$\tau\vert_{Z}$. Then $F_{i}\ne F_{j}$ if $i\ne j$, because $X$
satisfies the generality conditions and the cubic 3-fold $E$ is
smooth.

\begin{remark}
\label{remark:downstair-images-of-ramification-fibersd} The
surface $\pi\circ\alpha(D)=\Pi\subset\mathbb{P}^{4}$ is a
sufficiently general two-dimensional linear subspace passing
through the point $O$. The curve $\pi\circ\alpha(F_{i})\subset\Pi$
is one of three curves that are cut on the plane $\Pi$ by the
equation $f_{3}=0$. The line $\pi\circ\alpha(F_{i})$ is different
from the lines that are cut on $\Pi$ by the equation $f_{4}=0$.
Indeed, the plane $\Pi$ is sufficiently general, but the
polynomial $f_{4}$ is not divisible by the polynomial $f_{3}$ by
assumption. Therefore the fiber $F_{i}$ is smooth in the point of
intersection with the curve $Z$.
\end{remark}

The restriction morphism $\alpha\vert_{D}$ contracts the elliptic
curve $Z$ into the point $P$. The self-intersection of the curve
$Z$ on the surface $D$ is $-3$. The restriction
$\pi\vert_{\alpha(D)}$ is a cyclic triple cover of the plane $\Pi$
branched over a singular curve $\Pi\cap S$ of degree $6$ whose
singularities consist of the point, which is an ordinary triple
point on the curve $\Pi\cap S$.

Let $H\subset D$ be a curve that is cut on the surface $D$ by a
sufficiently general divisor in the linear system
$|\alpha^{*}(-K_{X})|$. The curve $H$ is smooth, the curve $H$ is
a three-section of the elliptic fibration $\tau$, the equality
$g(H)=4$ holds. Moreover, the curve $\pi\circ\alpha(H)\subset\Pi$
is a line. Let $C_{b}$ be a fiber of the elliptic fibration
$\tau:D\to\mathbb{P}^{1}$ over a point $b\in \mathbb{P}^{1}$. Then
$$
H^{2}=3, H\cdot Z=C_{b}^{2}=0, Z^{2}=-3, Z\cdot C_{b}=H\cdot
C_{b}=0
$$
on the surface $D$.

\begin{lemma}
\label{lemma:main-lemma-potential-density} For a very general
$\mathbb{C}$-point $b\in\mathbb{P}^{1}$ the equivalence
$$
3np-nH\vert_{C_{b}}\not\sim 0
$$
holds in $\mathrm{Pic}(C_{b})$ for every $n\in\mathbb{N}$, where
$p$ is one of the points of $Z\cap C_{b}$.
\end{lemma}

\begin{proof}
Let $T=Z\times_{\mathbb{P}^1}D$ be a fiber product and
$$
\chi:T\to D
$$
be an induced morphism. Then $\chi:T\to D$ is a cyclic triple
cover branched over the curves $F_{1}$, $F_{2}$, $F_{3}$. In
particular, the surface $T$ is singualr if and only if some fiber
$F_{i}$ is singular. However, the possible singularities of the
surface $T$ are easy to calculate in the case when we know the
type of the singular fiber $F_{i}$ of the elliptic fibration
$\tau$ (see \cite{BaPeVe84}). I

The surface $T$ is normal, and there is a well defined
intersection form of Weil divisors on the surface $T$ (see
\cite{Sak84}).

The fibration $\tau$ induces an elliptic fibration $\eta:T\to Z$
such that $\eta$ is a Jacobian fibration of the fibration $\tau$.
Indeed, the curve $\chi^{-1}(Z)$ splits into three irreducible
components, which are interchanged by the action of the group
$\mathbb{Z}_{3}$ on the surface $T$ that interchanges the fibers
of $\chi$. Let $\tilde{Z}$ be a component of the reducible curve
$\chi^{-1}(Z)$. Then $\tilde{Z}$ is a section of the fibration
$\eta$, and $\chi\vert_{\tilde Z}$ is an isomorphism.

Let $\tilde{H}=\chi^{-1}(H)$ and $L$ be a fiber of the fibration
$\eta$. Then the equalities
$$
\tilde{H}^{2}=9, \tilde{H}\cdot\tilde{Z}=L^{2}=0, \tilde{Z}\cdot L=1, \tilde{H}\cdot L=3%
$$
hold on the surface $T$. The curve $\tilde{Z}$ is smooth and
$\tilde{Z}\subset T\setminus\mathrm{Sing}(T)$, because the point
of intersection $F_{i}\cap Z$ is smooth on the fiber $F_{i}$ (see
Remark~\ref{remark:downstair-images-of-ramification-fibersd}).

The self-intersection $\tilde{Z}^{2}$ on $T$ can be calculated via
the adjunction formula, namely, we have $\tilde{Z}^{2}=-9$,
because $K_{T}\equiv 9L$.

It should be pointed out that in the case when the curve $Z$
passes through the singular points of the surface $T$ the
self-intersection $\tilde{Z}^{2}$ can be calculated using the
sub-adjunction formula with an appropriate different (see
\cite{Ko91}), which can be explicitly calculated for every type of
singular point.

For every $n\in\mathbb{N}$ we have
$$
3np-nH\vert_{C_{b}}\sim 0\iff (3n\tilde{Z}-n\tilde{H})\vert_{L_{a}}\sim 0\Rightarrow 3n\tilde{Z}-n\tilde{H}\equiv\Sigma, %
$$
where $C_{b}$ -- is a fiber of $\tau$ over a very general
$\mathbb{C}$-point $b\in\mathbb{P}^{1}$, $p$ is one of the
intersection points $Z\cap C_{b}$, $L_{a}$ is a fiber of $\eta$
over a very general $\mathbb{C}$-point $a\in Z$, and $\Sigma$ is a
divisor on the surface $T$ such that $\mathrm{Supp}(\Sigma)$ is a
union of the fibers of the elliptic fibration $\eta$.

Note, that all fibers of $\eta$ are irreducible, because all
fibers of $\tau$ are irreducible.

Suppose that the claim of the lemma is not true. Then the curves
$\tilde{Z}$, $\tilde{H}$, $L$ are linearly dependent in the group
$\mathrm{Div}(T)\otimes\mathbb{Q}\slash\equiv$. However, the
determinant of the matrix
$$
\bordermatrix{%
&                         &                        &                \cr %
& \tilde{Z}^{2}           &\tilde{H}\cdot\tilde{Z} & L\cdot\tilde{Z}\cr %
& \tilde{Z}\cdot\tilde{H} &\tilde{H}^{2}           & L\cdot\tilde{H}\cr %
& \tilde{Z}\cdot L        &\tilde{H}\cdot L        & L^{2}          \cr %
}=
\bordermatrix{%
&                          &                           &                          \cr %
& -9                       & 0                         & 1                        \cr %
& 0                        & 9                         & 3                        \cr %
& 1                        & 3                         & 0                        \cr %
}
$$
is $72\ne 0$, which contradicts to the linear dependence of the
curves $\tilde{Z}$, $\tilde{H}$, $L$.
\end{proof}

Now let us go from the surface $D$ back to the variety $U$. The
generality in the choice of the surface $D$ and
Lemma~\ref{lemma:main-lemma-potential-density} imply that
$$
3np+\alpha^{*}(nK_{X})\vert_{L_{p}}\not\sim 0
$$
in $\mathrm{Pic}(L_{p})$ for a very general $\mathbb{C}$-point
$p\in E$ and all $n\in\mathbb{N}$, where $L_{p}$ is a fiber of the
fibration $\psi:U\to\mathbb{P}^{3}$ over the point $p$.

For every $n\in\mathbb{N}$ let $\Phi_{n}\subseteq E$ be a subset
that is defined by the condition
$$
p\in \Phi_{n}\iff 3np\sim \alpha^{*}(-nK_{X})\vert_{L_{p}}
$$
in $\mathrm{Pic}(L_{p})$, where $L_{p}$ is a fiber of the elliptic
fibration $\psi$ over the point $\psi(p)$ such that the fiber
$L_{p}$ is smooth in a scheme-theoretic sense. Let
$\bar{\Phi}_{n}\subseteq E$ be a closure of the set $\Phi_{n}$ in
the Zariski topology. Then $\bar{\Phi}_{n}\ne E$ for every
$n\in\mathbb{N}$.

\begin{remark}
\label{remark:Merel-theorem} The set
$\Phi_{n}\setminus\cup_{i=1}^{n-1}\Phi_{i}$ does not contain
$\mathbb{F}$-points of the divisor $E$ for all natural numbers
$n>n(\mathbb{F})$ by Theorem~\ref{theorem:Merel}.
\end{remark}

The rational points are potentially dense on the divisor $E$ (see
Remark~\ref{remark:unirationality-of-cubic-threefold}). Thus we
can substitute $\mathbb{F}$ by its finite extension and assume
that $\mathbb{F}$-points of $E$ are Zariski dense.

Take an $\mathbb{F}$-point
$$
q\in
E\backslash\Big(\Delta\cup\cup_{i=1}^{n(\mathbb{F})}\bar{\Phi}_{i}\Big),
$$
where $\Delta$ is a Zariski closed subset of the divisor $E$
consisting of points that are contained in the singular fibers of
the elliptic fibration $\psi$. Let as before $L_{q}$ be a fiber of
$\psi$ over the point $\psi(q)$. Then $L_{q}$ and $\psi(q)$ are
defined over $\mathbb{F}$. Moreover, the curve $L_{q}$ is smooth.

By construction, the divisor $3q+\alpha^{*}(K_{X})\vert_{L_{q}}$
is defined over the field $\mathbb{F}$ and it is not a torsion in
$\mathrm{Pic}(L_{q})$. Therefore for any $n\in\mathbb{N}$ there is
a unique $\mathbb{F}$-point $q_{n}\in L_{b}$ such that
$$
q_{n}+(3n-1)q+\alpha^{*}(nK_{X})\vert_{L_{q}}\sim 0
$$
in $\mathrm{Pic}(L_{q})$ by the Riemann--Roch theorem.

We have $q_{i}\ne q_{j}$ if $i\ne j$. Hence the curve $L_{q}$ is
contained in the closure of all $\mathbb{F}$-points of the variety
$U$ in the Zariski topology for every $\mathbb{F}$-point $q$ in a
Zariski dense subset of the divisor $E$. Thus rational points are
dense on the varieties $U$ and $X$. It should be pointed out that
at certain point we substitute the field $\mathbb{F}$ by some its
finite extension in order to get the density of
$\mathbb{F}$-points on the divisor $E$. Hence
Proposition~\ref{proposition:explicit-potential-density} is
proved.

During the proof of
Proposition~\ref{proposition:explicit-potential-density} we
noticed that the surface $T$ is smooth if and only if each fiber
$F_{i}$ of the elliptic fibration $\tau$ is smooth. It is natural
to expect that the this is true for a sufficiently general $X$.
Indeed, the smoothness of the fiber $F_{i}$ is implied by the fact
that the line $\pi\circ\alpha(F_{i})$ intersects the ramification
hypersurface $S$ in three different points, one of which is the
point $O$. The latter condition can be easily expressed in terms
of the discriminant of the corresponding equation. Namely, it is
enough to require that the polynomials $f_{4}$ and
$f_{5}^{2}-4f_{4}f_{6}$ are not divisible by the irreducible
polynomial $f_{3}$.

Suppose that the polynomials $f_{4}$ and $f_{5}^{2}-4f_{4}f_{6}$
are not divisible by the irreducible polynomial $f_{3}$. Then the
divisor $E$ is a three-section of the elliptic fibration $\psi$
such that there is a smooth fiber $C$ of $\psi$ passing through
one of the ramification points of the restriction triple cover
$\psi\vert_{E}$. In the notations of the papers \cite{BoTsch98}
and \cite{BoTsch99} such multi-section is called saliently
ramified.

Let $C_{b}$ be a fiber of $\psi$ over a very general point
$b\in\mathbb{P}^{3}$, $p_{1}$ and $p_{2}$ be two different points
of  $C_{b}\cap E$. Then $p_{1}-p_{2}$ is not a torsion divisor on
$C_{b}$. Indeed, otherwise the torsion divisor $p_{1}-p_{2}$ goes
to a trivial divisor on $C$ when we $C_{b}\to C$. This arguments
can easily be put in algebraic form (see \cite{BoTsch98}). Now we
can prove the potential density of the rational points on $X$ in
the same way as in the proof of
Proposition~\ref{proposition:explicit-potential-density}, the only
difference is the following: we must generate $\mathbb{F}$-points
in the fibers of $\psi$ acting by the Jacobian fibration of $\psi$
without the usage the Riemann-Roch theorem (see \cite{BoTsch98}).

\end{document}